\documentclass[12pt]{article}
\usepackage{amsmath}
\usepackage{amsfonts}

\advance \textheight by \topskip
\advance \textheight by 1pt
\makeatother
\def\bea{\begin{eqnarray}}
\def\ena{\end{eqnarray}}

\def\lar{\longrightarrow}
\def\non{\nonumber}

\def\deg{\hbox{deg}}
\def\wt{\hbox{wt}}
\def\gr{\hbox{gr}}
\def\dim{\hbox{dim}}
\def\ch{\hbox{ch}}
\def\ker{\hbox{Ker}}

\def\sgn{\hbox{sgn}}
\def\ord{\hbox{ord}}

\def\Sym{\hbox{Sym}}

\def\tilp{\tilde{p}}

\def\tilp0{\tilde{p}_0}

\def\lar{\longrightarrow}

\def\ev{\hbox{ev}}
\def\evg{\hbox{ev}^{gr}}

\def\dum{du^{\hbox{max}}}

\newcommand{\qbc}[2]{
\left[
\begin{array}{c}{#1}\\{#2}\end{array}
\right]}

\newcommand{\qed}{\hbox{\rule[-2pt]{3pt}{6pt}}}
\newtheorem{prop}{Proposition}
\newtheorem{theorem}{Theorem}
\newtheorem{defn}{Definition}
\newtheorem{lemma}{Lemma}
\newtheorem{cor}{Corollary}
\newtheorem{conj}{Conjecture}
\newtheorem{remark}{Remark}

\title{
On Hyperelliptic Abelian Functions of Genus 3
}

\author{
Atsushi Nakayashiki\thanks{
e-mail: 6vertex@math.kyushu-u.ac.jp}\\
Department of Mathematics, Kyushu University\\
}

\date{}
\begin{document}
\maketitle

\vskip2mm
\centerline{
Dedicated to Mikio Sato on his eightieth birthday
}
\vskip15mm

\begin{abstract}
\noindent
The affine ring $A$ of the affine Jacobian variety $J(X)\backslash \Theta$
of a hyperelliptic curve of genus $3$ is studied as a ${\cal D}$ module.
A conjecture on the minimal ${\cal D}$-free resolution previously proposed
is proved in this case. As a by-product a linear basis of $A$ is explicitly
constructed in terms of derivatives of Klein's hyperelliptic $\wp$-functions.
\end{abstract}
\clearpage

\section{Introduction}
Let $X$ be a hyperelliptic curve of genus $3$, J(X) the Jacobian of $X$,
$\Theta$ the theta divisor, $A$ the affine ring of $J(X)\backslash\Theta$
and ${\cal D}$ the ring of holomorphic differential operators on $J(X)$.
The purpose of this paper is to determine the structure of $A$ as a 
${\cal D}$-module.

In general, for Jacobians, the non-linear differential equations satisfied
by elements of the affine ring are related with soliton equations 
\cite{B3,BEL1,M2}.
However the corresponding equations for non-Jacobians are not known.
The study of the ${\cal D}$-module structure of the affine ring $A$ of
an abelian variety is important from this point of view, 
since the results for Jacobians and non-Jacobians
can be compared in the same field.
 
In our previous paper \cite{NS1} a conjecture on the ${\cal D}$-free
resolution of $A$ is formulated in the case of hyperelliptic curves
of arbitrary genus. Up to now the conjecture is verified only for the cases
of genus 1 and 2. However those cases are contained in the generic case
where the theta divisor is non-singular.
In the case of a principally polarized abelian variety $(J,\Theta)$ with
$\Theta$ being non-singular, the ${\cal D}$-module structure of the affine ring is completely determined in \cite{CN}. Namely the minimal free resolution
is explicitly constructed. In the present case of genus $3$ the theta divisor
has an isolated singular point.
Thus it is the first case of the conjecture that is not contained 
in the generic one. 

Filtrations are important when one studies 
${\cal D}$-modules. We introduce two filtrations on $A$,
pole and KP-filtrations. The pole filtration is
defined by the order of poles on $\Theta$ and it can be defined
for other abelian varieties than Jacobians.
To define KP-filtration we use Klein's hyperelliptic sigma function
$\sigma(u)=\sigma(u_1,...,u_g)$ \cite{K1,K2,BEL1,N2}. We consider $\Theta$ 
as the zero set of $\sigma(u)$. Let 
\bea
&&
\wp_{i_1...i_n}(u)=-\partial_{i_1}\cdots\partial_{i_n}\log\,\sigma(u),
\qquad
\partial_i=\frac{\partial}{\partial u_i}.
\non
\ena
For $n\geq 2$ $\wp_{i_1...i_n}(u)$ is contained in $A$ and conversely 
$A$ is generated by $\{\wp_{i_1...i_n}(u)\}$ as a ring. Assign
degree $\sum_{j=1}^{n}(2i_j-1)$ to $\wp_{i_1...i_n}(u)$. Then the 
KP-filtration $\{A_n\}$ is defined by specifying $A_n$ to be the
vector space generated by elements of degree at most $n$. The KP-filtration
is specific to Jacobians and is related to integrable systems known
as KP-hierarchy \cite{DJKM,Sat}.

The KP-filtration seems to be a proper filtration to study
the affine ring of a Jacobian. In fact the result on the
character \cite{NS1,NS2}, which is the generating function of 
the dimensions of homogeneous components of the associated graded 
ring, with respect
to the KP-filtration manifests a remarkable consistency with
other results and constructions, such as 
the results on the cohomology groups of affine Jacobians \cite{N1}.

Nevertheless we use the pole filtration for the proof of the conjecture
in the present case. There are three reasons for this. 
One is that the pole filtration
can be localized and the sheaf cohomology arguments can be applied. In fact
the local structure is inherited to the global structure and it plays
a decisive role to determine the ${\cal D}$-module structure of $A$.
The second is that we are interested in describing explicitly 
 a basis of abelian functions with poles of order at most $n$. 
This is for the sake of the application
to finding
explicit relations among abelian functions, such as generalizations
of Frobenius-Stickelberger's formula \cite{B2,BEL1,EEMOP,EEP}. 
The third is 
the lack of the technical device to treat the KP-filtration. For example
one can not define the corresponding filtration locally on $J(X)$.
It is important to develop intrinsic geometric understanding
of the KP filtration for the further study. 

Let $A=\cup A(n)$  and $A=\cup A_n$ be the pole and the KP filtrations.
Denote by $\gr^P\,A$ and $\gr^{KP}\,A$ the graded rings associated with 
the pole and KP filtrations respectively. They also become ${\cal D}$-modules.
We prove that $\gr^P\,A$ is not finitely generated over ${\cal D}$
and analyze how it is not finitely generated. It is shown that
the elements of $A(n)$ is not contained in ${\cal D}^{(1)}A(n-1)$
but is contained in ${\cal D}^{(2)}A(n-1)$, where
${\cal D}^{(k)}$ is the space of differential operators of order at most $k$.
Namely some elements of $A(n)/A(n-1)$ are not obtained by differentiating once 
the elements of $A(n-1)$ but are obtained by 
differentiating twice and taking linear combinations.
This phenomenon is a result of the existence
of the singularity of the theta divisor. 
It gives us an insight, for more general cases, on what happens
and what we should prove if $\Theta$ is singular.
To establish such results we need
to study the residue sheaves supported on the singular locus of 
the theta divisor \cite{AH,Z}. 
To this end 
we use Taylor expansion of the sigma function. In fact one of the
important properties of the sigma function is that the series expansion
is known explicitly \cite{B2,BL,BEL1,N2}.
The first term of the expansion is given by Schur function 
corresponding to the partition 
determined from the gap sequence at $\infty$ of $X$. 
In the present case the partition is $(3,2,1)$ and the 
corresponding Schur function is
\bea
&&
S(u)=u_1u_3-u_2^2-\frac{1}{3}u_1^3u_2+\frac{1}{45}u_1^6.
\non
\ena
The zero set of $S(u)$ has a simple singularity of type $A_1$ at the origin.
It implies, in particular, that $\sigma(u)$ is transformed to some 
canonical polynomial defining $A_1$-type singularity 
around the origin by taking a suitable local coordinate system. 
With the aid of the explicit form of 
the local defining equation we can analyze residue sheaves in detail.
Then the differential property of $A(n)$ mentioned above can be proved
by taking cohomology.

We can deduce from the results on $\gr^P\,A$ that 
$A$ is finitely generated over ${\cal D}$ although $\gr^P\,A$ is not. 
Moreover generators can be taken as elements in $A$ of the form
\bea
&&
1,
\quad
\wp_{ij}(u),
\quad
\left|
\begin{array}{cc}
\wp_{i_1j_1}(u)&\wp_{i_1j_2}(u)\\
\wp_{i_2j_1}(u)&\wp_{i_2j_2}(u)\\
\end{array}
\right|,
\quad
\left|
\begin{array}{ccc}
\wp_{11}(u)&\wp_{12}(u)&\wp_{13}(u)\\
\wp_{21}(u)&\wp_{22}(u)&\wp_{23}(u)\\
\wp_{31}(u)&\wp_{32}(u)&\wp_{33}(u)\\
\end{array}
\right|.
\label{co-basis}
\ena
 Next we derive ${\cal D}$-linear relations among 
derivatives of (\ref{co-basis}) 
and determine a linear basis of $A$. 
With the help of this linear basis $\gr^{KP}\,A$ is proved to be
generated by (\ref{co-basis}) over ${\cal D}$.
Once this is established the conjecture on ${\cal D}$-free resolutions
 of $\gr^{KP}\,A$ and $A$ are proved to be true. 
In this way we determine
the ${\cal D}$-module structure and a ${\mathbb C}$-basis of $A$. 

The present paper is organized in the following manner.
In section 2 we review the definition and fundamental properties
of the hyperelliptic sigma function following \cite{BEL1,N2}.
The matrix construction of the affine hyperelliptic Jacobian is reviewed
and the KP-filtration is introduced in section 3. In section 4
the conjecture of \cite{NS1} is reviewed and the main result of this
paper is given. The local differential structure of sheaves is studied by 
analyzing the local defining equation of the theta divisor near the singular
point in section 5. In section 6 cohomology groups of sheaves
with higher order poles are studied. 
It is shown that $A$ is finitely generated over ${\cal D}$
while $\gr^P\,A$ is not.
The explicit description of the cohomology group 
$H^3(J(X)\backslash\Theta,{\mathbb C})$ is reviewed in section 7.
In section 8 the addition theorem of the genus 3 hyperelliptic sigma
function due to H. F. Baker is reviewed and a basis of $A(2)$ is determined
in terms of the cohomology of $J(X)\backslash\Theta$ given 
in the previous section. Bases of Abelian functions of lower order poles are
studied in section 9 and 10. As a consequence it is shown that $A$ is generated
by representatives of the cohomology group 
$H^3(J(X)\backslash\Theta,{\mathbb C})$ given in section 7.
In section 11 a linear basis of $A$ is determined as a subset of derivatives
of the generators given in section 10. Finally a basis of $\gr^{KP}\,A$ is 
determined and the proof of the conjecture is given in section 12. In section
13 some remarks are given and remaining problems are discussed.

\section{Sigma Function}
In this section we recall the definition and fundamental properties of the
hyperelliptic sigma function \cite{K1,K2}. See \cite{BEL1,N2} for more details.

Consider the hyperelliptic curve defined by the equation
\bea
&&
y^2=f(x),
\qquad
f(x)=\sum_{i=0}^{2g+1}\lambda_i x^i,
\quad
\lambda_{2g+1}=4.
\non
\ena
We assume that $f(x)$ has no multiple roots.
Let $X$ be the corresponding compact Riemann surface of genus $g$ and
\bea
&&
du_i=\frac{x^{g-i}dx}{y},
\quad
i=1,...,g
\non
\ena
a basis of holomorphic one forms on $X$. 
We consider the second kind differentials defined by
\bea
&&
dr_i=
\sum_{k=g+1-i}^{g+i}(k-g+i)\lambda_{k+g+2-i}
\frac{x^kdx}{4y},
\quad
i=1,...,g.
\non
\ena
Being considered as elements of $H^1(X,{\mathbb C})$, $\{du_i,dr_i\}$ 
forms a symplectic basis with respect to the intersection form $\circ$:
\bea
&&
du_i\circ du_j=dr_i\circ dr_j=0,
\quad
du_i\circ dr_j=\delta_{ij}.
\label{co-symp}
\ena
By specifying a symplectic basis of the homology group $H_1(X,{\mathbb Z})$
we define period matrices:
\bea
&&
2\omega_1=\left(\int_{\alpha_j}du_{i}\right),
\quad
2\omega_2=\left(\int_{\beta_j}du_{i}\right),
\quad
-2\eta_1=\left(\int_{\alpha_j}dr_{i}\right),
\quad
-2\eta_2=\left(\int_{\beta_j}dr_{i}\right),
\non
\ena
and $\tau=\omega_1^{-1}\omega_2$.

Let $p_n(T)$ be the polynomial of $\{T_i\}$
defined by
\bea
&&
\exp(\sum_{n=1}^\infty T_nk^n)=\sum_{n=0}^\infty p_n(T)k^n.
\non
\ena
For a partition $\lambda=(\lambda_1,...,\lambda_l)$ 
define Schur function $S_\lambda(T)$ by
\bea
&&
S_\lambda(T)=\det (p_{\lambda_i-i+j}(T))_{1\leq i,j\leq l}.
\non
\ena
\vskip5mm

\noindent
{\bf Example}
\hskip5mm
$\displaystyle{S_{(2,1)}(T)=-T_3+\frac{T_1^3}{3}}$,
\hskip5mm
$\displaystyle{S_{(3,2,1)}(T)=T_1T_5-T_3^2-\frac{1}{3}T_1^3T_3
+\frac{1}{45}T_1^6}$.
\vskip10mm

We assign degree $-i$ to $T_i$.
Then $S_\lambda(T)$ is homogeneous of degree $-|\lambda|$, 
where $|\lambda|=\lambda_1+\cdots+\lambda_l$.
For each $g\geq 1$ we define the partition
\bea
&&
\lambda(2,2g+1)=(g,g-1,...,1).
\non
\ena
The function $S_{\lambda(2,2g+1)}(T)$
becomes a polynomial of $T_1,T_3,...,T_{2g-1}$. 
 Consider the variables $u_i$, $1\leq i\leq g$ and
assign the degree as $\deg\,u_i=-(2i-1)$.

Let $\delta'+\tau\delta''$ with $\delta',\delta''\in 1/2{\mathbb Z}^g$ 
be the Riemann constant with respect to the base point $\infty$.

\begin{defn}\label{def-1}
The fundamental sigma function or simply the sigma function
$\sigma(u)$ is the holomorphic function on ${\mathbb C}^g$ of the
variables $u={}^t(u_1,...,u_g)$ which satisfies the following
conditions.
\vskip2mm
\noindent
(i) For any $m_1, m_2 \in {\mathbb Z}^g$,
\vskip2mm
\hskip3mm
$
\displaystyle{
\sigma(u+2\omega_1m_1+2\omega_2m_2)
=(-1)^{{}^tm_1m_2+2({}^t\delta'm_1-{}^t\delta''m_2)}
}
$
\bea
&&
\qquad\qquad\qquad\qquad
\times
\exp\left(
{}^t(2\eta_1m_1+2\eta_2m_2)(u+\omega_1m_1+\omega_2m_2)
\right)\sigma(u).
\non
\ena

\vskip2mm
\noindent
(ii) The expansion of $\sigma(u)$ at the origin is of the form
\bea
&&
\sigma(u)=S_{\lambda(2,2g+1)}(T)|_{T_{2i-1}=u_i}+\sum_{d} f_d(u),
\label{normalization}
\ena
where $f_d(u)$ is a homogeneous polynomial of degree $d$ and the sum
is over integers $d$ satisfying $d<-|\lambda(2,2g+1)|$.
\end{defn}

The sigma function can be written in terms of Riemann's theta function
as
\bea
&&
\sigma(u)=
C\exp\left(\frac{1}{2}{}^tu\eta_1\omega_1^{-1}u\right)
\theta\qbc{\delta'}{\delta''}((2\omega_1)^{-1}u,\tau),
\label{sigma-theta}
\ena
where $C$ is some constant specified by (ii) of Definition \ref{def-1}
(see \cite{BEL1} for the explicit formula for $C$). 

Let
\bea
&&
J(X)={\mathbb C}^g/2\omega_1{\mathbb Z}^g+2\omega_2{\mathbb Z}^g
\non
\ena
be the Jacobian variety of $X$ and $\Theta$ the divisor defined by the 
zero set of $\sigma(u)$. 
We call $\Theta$ the theta divisor throughout this paper.
\vskip10mm

\section{Affine Jacobian}
In this section we review the matrix construction 
of the affine Jacobian variety $J(X)\backslash\Theta$ 
due to Jacobi and Mumford \cite{M2,NS1} 
and give a description of the generators of the affine 
ring in terms of the sigma function \cite{BEL1,M2}.
We mainly follow the notation in \cite{NS1}.

Let ${\cal L}$ be the set of matrices of the form
\bea
&&
L(x)=
\left[
\begin{array}{cc}
a(x) & b(x)\\
c(x) & -a(x)\\
\end{array}
\right],
\non
\\
&&
a(x)=\sum_{i=1}^ga_{2i+1}x^{g-i},
\quad
b(x)=\sum_{i=0}^g b_{2i}x^{g-i},
\quad
c(x)=\sum_{i=0}^{g+1} c_{2i}x^{g+1-i},
\quad
b_0=1, c_0=4.
\non
\ena
We set $a_1=0$. 
Here the choice of $b_0$ and $c_0$ corresponds to the choice of 
the coefficient of the highest degree term of $f(x)$ below.
We identify ${\cal L}$ with the affine space
${\mathbb C}^{3g+1}$ by the map
\bea
&&
L(x)\mapsto
(a_3,...,a_{2g+1},b_2,...,b_{2g},c_2,...,c_{2g+2}).
\non
\ena
For a polynomial $f(x)=\sum_{i=0}^{2g+1}\lambda_i x^i$, $\lambda_{2g+1}=4$,
consider the equation
\bea
-\det\,L(x)=f(x).
\label{eq-n21}
\ena
It gives a set of equations for $a_i,b_j,c_k$.
Let ${\cal L}_f$ be the set of elements of ${\cal L}$ satisfying
(\ref{eq-n21}). 

\begin{theorem}\label{th-n21}{\rm \cite{M2}}
If $f(x)$ does not have multiple roots, ${\cal L}_f$ is an affine algebraic
variety and is isomorphic to $J(X)\backslash \Theta$.
\end{theorem}

Let 
\bea
&&
{\bf A}={\mathbb C}[a_{2i+1},b_{2j},c_{2k}\,|\, 1\leq i\leq g,
 1\leq j\leq g, 1\leq k\leq g+1]
\non
\ena
be the polynomial ring of $3g+1$ variables and $I_f$ the ideal generated
by the coefficients of (\ref{eq-n21}). Then 
\bea
&&
A_f={\bf A}/I_f
\label{eq-n22}
\ena
is the affine ring of ${\cal L}_f$. Notice that $A_f$ is generated by
$a_i,b_j$ since $c_k$ is expressed as a polynomial of $a_i,b_j$ by 
(\ref{eq-n21}) \cite{NS1}.

The affine ring of $J(X)\backslash\Theta$ is isomorphic to the 
ring of meromorphic functions on $J(X)$ with poles only on $\Theta$.
Meromorphic functions on $J(X)$ are identified with those on ${\mathbb C}^g$
that are periodic with respect to the lattice 
$2\omega_1{\mathbb Z}^g+2\omega_2{\mathbb Z}^g$.
Such periodic functions can be constructed as logarithmic derivatives of the
sigma function for example. Let
\bea
&&
\wp_{i_1...i_n}(u)=-\partial_{i_1}\cdots\partial_{i_n}\log\,\sigma(u),
\quad
\partial_i=\frac{\partial}{\partial u_i}.
\non
\ena
For $n\geq 2$ $\wp_{i_1...i_n}(u)$ becomes an element of the affine ring
of $J(X)\backslash\Theta$. According to Theorem \ref{th-n21} $a_i,b_j$ 
should be described as meromorphic functions on $J(X)$.
The result is known as \cite{BEL1}
\bea
&&
b_{2i}=-\wp_{1i}(u),
\qquad
a_{2j+1}=\wp_{11j}(u).
\label{eq-n22-1}
\ena

In the following we fix $f(x)$ and denote $A_f$ simply by $A$.
We introduce a filtration on $A$ using the relation (\ref{eq-n22}).
Define a grading on ${\bf A}$ by
\bea
&&
\deg\,a_i=i,
\qquad
\deg\,b_i=i,
\qquad
\deg\,c_i=i.
\non
\ena
Let
\bea
&&
{\bf A}=\oplus_{n=0}^\infty {\bf A}_n,
\qquad
{\bf A}_0={\mathbb C},
\non
\ena
be the homogeneous decomposition of ${\bf A}$ and $\pi:{\bf A}\rightarrow A$
the projection.
We set
\bea
&&
A_n=\pi(\oplus_{d=0}^n {\bf A}_d)
\non
\ena
for $n\geq 0$ and $A_n=0$ for $n<0$. Obviously $\{A_n\}$ defines
an increasing filtration of $A=\cup_{n=0}^\infty A_n$ which we call
KP-filtration. Let $\gr^{KP}\,A$ be the associated graded ring
\bea
&&
\gr^{KP}\,A=\oplus_{n=0}^\infty \gr^{KP}_n\,A,
\quad
\gr^{KP}_n\,A=A_n/A_{n-1}.
\non
\ena

\begin{lemma}\label{lem-n22}
For $n\geq 2$ and $i_1,...,i_n\in \{1,...,g\}$ we have
\bea
&&
\wp_{i_1...i_n}\in A_{\sum_{j=1}^n(2i_j-1)}.
\non
\ena
\end{lemma}

To prove this lemma we first describe the action of $\partial_i$
on $a_j$, $b_k$. The translation invariant vector field $D_i$ on $J(X)$ is
constructed in \cite{M2}. It gives
\bea
&&
D_l(a_{2k+1})=\frac{1}{4}\sum(b_{2i}c_{2j+2}-b_{2j}c_{2i+2})-b_{2k}b_{2l},
\non
\\
&&
D_l(b_{2k})=\frac{1}{2}\sum(a_{2i+1}b_{2j}-a_{2j+1}b_{2i}),
\non
\\
&&
D_l(c_{2k+2})=
\frac{1}{2}\sum(c_{2i+2}a_{2j+1}-c_{2j+2}a_{2i+1})+2b_{2l}a_{2k+1},
\label{Dl}
\ena
where all sums are over $(i,j)$ satisfying
\bea
&&
i+j=k+l-1,
\qquad
i\geq\max(k,l),
\qquad
j\leq \min(k,l)-1.
\non
\ena
Notice that, due to the coefficient $4$ of $x^{2g+1}$ in $f(x)$,
 the coefficients in the right hand side of (\ref{Dl})
 are different from those in \cite{M2}. We have
\bea
&&
D_l(b_2)=\frac{1}{2}a_{2l+1}.
\non
\ena
In terms of $\wp_{i_1...i_n}$
\bea
&&
D_l\wp_{11}(u)=-\frac{1}{2}\wp_{11l}(u).
\label{eq-n23}
\ena

\begin{lemma}\label{lem-n23}
$\displaystyle{D_l=-\frac{1}{2}\partial_l}$.
\end{lemma}
\vskip2mm
\noindent
{\it Proof.} The equation (\ref{eq-n23}) is written as
\bea
&&
(D_l+\frac{1}{2}\partial_l)(\wp_{11}(u))=0.
\non
\ena
Therefore it is sufficient to prove the following statement:
if an invariant vector field $D=\sum_{i=1}^g\alpha_i \partial_i$ satisfies
\bea
&&
D\wp_{11}(u)=0,
\label{eq-n24}
\ena
then $D=0$.
In fact (\ref{eq-n24}) implies 
\bea
&&
\frac{2\sigma_1^2D\sigma}{\sigma^3}-
\frac{D(\sigma_1^2)+\sigma_{11}D\sigma}{\sigma^2}
+\frac{D(\sigma_{11})}{\sigma}=0,
\non
\ena
where $\sigma_1=\partial_1\sigma$, $\sigma_{11}=\partial_1^2\sigma$.
It means that $\sigma_1^2D\sigma/\sigma$ is holomorphic.

By claim (i) of Lemma \ref{ord-zeta}
\bea
&&
\wp_{11}(u)=-\frac{\sigma_1^2}{\sigma^2}+\frac{\sigma_{11}}{\sigma}
\non
\ena
has poles of order two on $\Theta$. 
Since $\Theta$ is irreducible, $D\sigma/\sigma$ is holomorphic.
Then
\bea
\partial_1\frac{D\sigma}{\sigma}=-\sum_{i=1}^g\alpha_i\wp_{1i}(u)
\non
\ena
is holomorphic on $J(X)$. Thus it is a constant. Since $\{1,\wp_{ij}\,|1\leq i\leq j\leq g\,\}$ is linearly independent by Proposition \ref{cohomology-basis}, $\alpha_i=0$ for any $i$. Thus Lemma \ref{lem-n23}
is proved.
\qed

\vskip5mm
\noindent
{\it Proof of Lemma \ref{lem-n22}}
\vskip2mm
\noindent
Since $b_2=-\wp_{11}(u)$, we have
\bea
&&
\wp_{11}(u)\in A_2.
\non
\ena
The formulae (\ref{Dl}) shows that
\bea
&&
D_lA_n\subset A_{n+2l-1}.
\label{eq-n25}
\ena
The relation of $\sigma(u)$ to the $\tau$-function of the KP-hierarchy \cite{EG,N3}, which in fact reduces to the KdV hierarchy in the present case,
  implies that $\wp_{ij}(u)$ is expressed as a homogeneous 
polynomial of $\wp_{11}$, $\wp_{111}$,... of degree $(2i-1)+(2j-1)$ 
modulo $A_{2(i+j)-3}$, where the homogeneity is with respect to
the degree $\deg\,\partial_1^i\wp_{11}=i+2$. 
Thus, by (\ref{eq-n25}) and the obvious relation $A_mA_n\subset A_{m+n}$, 
we have
\bea
&&
\wp_{i_1i_2}(u)\in A_{(2i_1-1)+(2i_2-1)}.
\non
\ena
Applying $\partial_{i_3}\cdots\partial_{i_n}=(-2)^{n-2}D_{i_3}\cdots D_{i_n}$
to $\wp_{i_1i_2}(u)$ we get the desired result. \qed

In general, for a graded vector space $S=\oplus_n S_n$, we define
the character of $S$ as the generating function of the dimensions of 
homogeneous components:
\bea
&&
\ch(S)=\sum q^n\dim\, S_n.
\non
\ena
To give a formulae for the character of $\gr^{KP}\,A$ we introduce
the notation:
\bea
&&
[n]_p=1-p^n,
\qquad
[n]_p!=\prod_{i=1}^n\,\,[\,i\,]_p,
\qquad
[n+\frac{1}{2}]_p!=\prod_{i=0}^n\,\,[\,i+\frac{1}{2}\,\,]_p,
\non
\ena
for a non-negative integer $n$.

\begin{theorem}\label{th-n24}{\rm \cite{NS1}}
The following formula is valid:
\bea
&&
\ch(\gr^{KP}\,A)=
\frac{[\,\frac{1}{2}\,]_{q^2}[2g+1]_{q^2}!}
{[g]_{q^2}![g+1]_{q^2}![g+\frac{1}{2}]_{q^2}!}.
\non
\ena
\end{theorem}
 
In this paper we also consider another filtration on $A$, the pole filtration,
defined as follows. For $a\in A$ we denote by $\ord\,a$ the order of poles on 
$\Theta$. Set
\bea
&&
A(n)=\{\,a\in A\,|\, \ord\, a\leq n\,\}.
\label{eq-n26}
\ena
Then $\{A(n)\}$ defines an increasing filtration on $A$.
Notice that $A(0)=A(1)={\mathbb C}$. The graded ring associated
with this filtration is denoted by $\gr^{P}\,A$:
\bea
&&
\gr^P\,A=\oplus_{n=0}^\infty\gr^P_n\,A,
\qquad
\gr^P_n\,A=A(n)/A(n-1).
\non
\ena
It is obvious that the following relation holds:
\bea
&&
\partial_i A(n)\subset A(n+1).
\label{der-pole}
\ena

\section{Abelian Functions as a ${\cal D}$-module}
Let ${\cal D}={\mathbb C}[\partial_1,...,\partial_g]$ be the ring of holomorphic differential operators on $J(X)$. 
As observed in the previous section the affine ring $A$ of 
$J(X)\backslash\Theta$ becomes a ${\cal D}$-module.
The relations (\ref{eq-n25}) and (\ref{der-pole}) imply that 
$\gr^{KP}\,A$ and $\gr^P\,A$ become also  
${\cal D}$-modules. In this section we recall the conjecture on the 
${\cal D}$-module structure of $A$ and $\gr^{KP}\,A$ proposed in \cite{NS1}.

Let 
\bea
&&
V=\oplus_{i=1}^g {\mathbb C}\epsilon_i \oplus \oplus_{i=1}^g{\mathbb C}\mu_i
\non
\ena
be the vector space of dimension $2g$ with the basis $\{\epsilon_i,\mu_i\}$.
Consider the two form
\bea
&&
\omega=\sum_{i=1}^g \epsilon_i\wedge \mu_i\in \wedge^2 V,
\non
\ena
and set
\bea
&&
W^k=\frac{\wedge^k V}{\omega\wedge^{k-2} V}
\quad
k\geq 2,
\quad
W^1=V,
\quad
W^0={\mathbb C}.
\label{wk}
\ena
We define a grading on $V$ by assigning
\bea
&&
\deg\, \epsilon_i=-(2i-1),
\qquad
\deg\,\mu_i=2i-1.
\non
\ena
Then $\wedge^k V$ for $k\geq 2$ is naturally graded and $\deg\,\omega=0$.
Thus $W^k$ is also graded as the quotient of two graded spaces.

Define the map
\bea
&&
d:{\cal D}\otimes \wedge^k V \lar {\cal D}\otimes \wedge^{k+1} V
\non
\ena
by
\bea
&&
d(P\otimes \nu)=\sum \partial_i P\otimes (\epsilon_i\wedge \nu),
\quad
P\in {\cal D},
\quad
\nu\in \wedge^k V,
\non
\ena
and the map
\bea
&&
\omega:
{\cal D}\otimes \wedge^k V \lar {\cal D}\otimes \wedge^{k+2} V
\non
\ena
by
\bea
&&
\omega(P\otimes \nu)=P\otimes (\omega\wedge \nu).
\non
\ena
We specify a grading on ${\cal D}$ by
\bea
&&
\deg\, \partial_i=2i-1.
\non
\ena
The space ${\cal D}\otimes \wedge^k V$ naturally inherits a grading
from ${\cal D}$ and $\wedge^k V$. The maps $d$ and $\omega$ preserve
the grading since $\deg\,\omega=0$ and $\deg\,d=\deg\,\sum \partial_i\otimes \epsilon_i=0$. Obviously $d$ and $\omega$ commute and $d^2=0$.
Therefore $d$ induces a map
\bea
&&
d:{\cal D}\otimes W^k \lar {\cal D}\otimes W^{k+1},
\non
\ena
and defines the complex $({\cal D}\otimes W^{\bullet},d)$:
\bea
&&
0\lar {\cal D}\lar {\cal D}\otimes W^1 \lar\cdots\lar {\cal D}\otimes W^g\lar 0.\non
\ena

\begin{prop}{\rm \cite{NS1}}
The complex $({\cal D}\otimes W^{\bullet},d)$ is exact at 
${\cal D}\otimes W^{k}$, $k\neq g$.
\end{prop}

Let
\bea
&&
T^\ast=\sum_{i=1}^g {\mathbb C}du_i
\non
\ena
be the space of holomorphic one forms on $J(X)$.
We define the map
\bea
&&
\ev:{\cal D}\otimes W^g \lar A\otimes \wedge^g T^\ast,
\non
\ena
as follows.
Let
\bea
&&
\zeta_i(u)=\partial_i\log\,\sigma(u),
\qquad
\zeta_{ij}(u)=-\wp_{ij}(u)=\partial_i\partial_j\log\,\sigma(u).
\non
\ena
Then
\bea
&&
d\zeta_i=\sum_{j=1}^g\zeta_{ij}(u)du_j\in
A\otimes T^\ast.
\non
\ena
We set
\bea
&&
\dum=du_1\wedge\cdots\wedge du_g.
\non
\ena
Since $\wedge^gT^\ast={\mathbb C}\dum$, $A\otimes \wedge^gT^\ast$ becomes
a ${\cal D}$-module by
\bea
&&
P(a\otimes \dum)=P(a)\otimes \dum.
\non
\ena
As a ${\cal D}$-module $A\otimes \wedge^gT^\ast$ and $A$ are isomorphic.
For $I=(i_1,...,i_r)\in\{1,...,g\}^r$ we use the notation like
\bea
&&
\epsilon_I=\epsilon_{i_1}\wedge\cdots\wedge\epsilon_{i_r}.
\non
\ena
We define
\bea
&&
\ev
\left(
P \otimes 
(
\mu_I
\wedge
\epsilon_J
)
\right)
=P\left(
d\zeta_I
\wedge 
du_J
\right),
\non
\ena
where $P\in {\cal D}$,
$I=(i_1,...,i_r)\in \{1,...,g\}^r$ and 
$J=(j_{r+1},...,j_g)\in \{1,...,g\}^{g-r}$. 

The map $\ev$ can be written explicitly in terms of $\zeta_{ij}$.
Let $J^c=(j_1,...,j_r)$ be defined such that $j_1<\cdots<j_r$ and
$\{1,...,g\}\backslash J=\{j_1,...,j_r\}$.
\bea
&&
d\zeta_I\wedge du_J=sgn(J^c,J)(I;J^c)\dum,
\non
\\
&&
(I;J^c)=\det(\zeta_{i_kj_l})_{1\leq k,l\leq r},
\non
\ena
where $sgn(J^c,J)$ is the sign of the permutation $(J^c,J)$.
Notice that $(i;j)=\zeta_{ij}$.
Then we have
\bea
&&
\ev\left(
P\otimes (\mu_I\wedge \epsilon_{J})
\right)
=
\sgn(J^c,J)P\left((I;J^c)\right)\dum.
\non
\ena
We also define the graded version $\evg$ of the map $\ev$.
To this end let us define a grading on $T^\ast$ by
\bea
&&
\deg\, du_i=-(2i-1).
\non
\ena
Then $\gr^{KP}\,A\otimes \wedge^g T^\ast$ is graded.
Let
\bea
&&
{\cal D}\otimes W^g=\oplus \left({\cal D}\otimes W^g\right)_n,
\quad
\gr^{KP}\,A\otimes \wedge^g T^\ast=
\oplus \left(\gr^{KP}\,A\otimes \wedge^g T^\ast\right)_n,
\non
\ena
be the homogeneous decompositions. Notice that
\bea
&&
\left(\gr^{KP}\,A\otimes \wedge^g T^\ast\right)_n=
\gr^{KP}_{n+g^2}\,A\otimes \dum.
\non
\ena
We have
\bea
&&
\deg(\mu_I\wedge \epsilon_J)=\sum_{k=1}^r(2i_k-1)-\sum_{k=r+1}^g(2j_k-1)
=:d_{I,J}.
\non
\ena
On the other hand a calculation shows that
\bea
&&
(I;J^c)\in A_{d_{I,J}+g^2}.
\non
\ena
Thus
\bea
&&
\ev\left(({\cal D}\otimes W^g)_n\right)
\subset
A_{n+g^2}\otimes \dum.
\non
\ena
Composing $\ev$ and the projection $A_{n+g^2}\rightarrow \gr^{KP}_{n+g^2}\,A$ 
one can define
\bea
&&
\evg_n:({\cal D}\otimes W^g)_n \lar 
\left(\gr^{KP}\,A\otimes \wedge^g T^\ast\right)_n.
\non
\ena
We set
\bea
&&
\evg=\oplus_n \evg_n:{\cal D}\otimes W^g \lar 
\gr^{KP}\,A\otimes \wedge^g T^\ast.
\non
\ena

\begin{conj}\label{conj1}{\rm \cite{NS1}}
The map $\evg$ is surjective. In other words $\gr^{KP}\,A$ is generated, as
a ${\cal D}$-module, by $1\in \gr^{KP}_{0}\,A$ and 
$(I;J^c)\in \gr^{KP}_{d_{I,J}+g^2}\,A$, $r\geq1$, 
$I=(i_1,...,i_r)\in\{1,...,g\}^r$, 
$J=(i_{r+1},...,i_g)\in\{1,...,g\}^{g-r}$.
\end{conj}

\begin{cor}{\rm \cite{NS1}}
If Conjecture \ref{conj1} is true, the following two complexes are exact
and give ${\cal D}$-free resolutions of $\gr^{KP}\,A$ and $A$ respectively:
\bea
&&
0\lar {\cal D}
\stackrel{d}{\lar} 
{\cal D}\otimes W^1
\stackrel{d}{\lar} 
\cdots
\stackrel{d}{\lar} 
{\cal D}\otimes W^g
\stackrel{\evg}{\lar} 
\gr^{KP}\,A\otimes \wedge^g T^\ast\lar 0,
\non
\\
&&
0\lar {\cal D}
\stackrel{d}{\lar} 
{\cal D}\otimes W^1
\stackrel{d}{\lar} 
\cdots
\stackrel{d}{\lar} 
{\cal D}\otimes W^g
\stackrel{\ev}{\lar} 
A\otimes \wedge^g T^\ast\lar 0.
\non
\ena
\end{cor}

For $g=1$ the conjecture is obvious since $\{1,\wp(u),\wp'(u),...\}$
is a linear basis of $A$. Here $\wp(u)=\wp_{11}(u)$ is  
Weierstrass' elliptic function. 
For $g=2$ the conjecture follows from Example 9.2 in \cite{CN}. 

In this paper we prove

\begin{theorem}\label{main-th}
Conjecture \ref{conj1} is true for $g=3$.
\end{theorem}

\section{Local Structure}
From this section until the end of the paper we assume $g=3$. 
In this case the only singularity of the divisor $\Theta$ is the point 
corresponding to $(u_1,u_2,u_3)=(0,0,0)$, which we denote $0\in J(X)$.

For $p\geq 0$, $n\in {\mathbb Z}$, let $\Omega^p$ be the sheaf of germs 
of holomorphic $p$-forms on $J(X)$, $\Omega^p(n)$ $(n\geq 0)$ the sheaf
 of germs of meromorphic $p$-forms on $J(X)$ which have poles only on $\Theta$
of order at most $n$, $\Omega^p(n)$ $(n<0)$ the sheaf of germs of holomorphic
$p$-forms on $J(X)$ which have zeros on $\Theta$ of order at least $-n$.
We set ${\cal O}=\Omega^0$, ${\cal O}(n)=\Omega^0(n)$ and 
$\gr_n\,\Omega^p=\Omega^p(n)/\Omega^p(n-1)$. 
Since $\Omega^p$ is a free ${\cal O}$-module, $\gr_n\,\Omega^p\simeq 
\gr_n\,{\cal O}\otimes \Omega^p$.

The exterior differentiation defines a map 
$d: \Omega^p(n)\lar \Omega^{p+1}(n+1)$.
It induces a map $d: \gr_n\,\Omega^p\lar \gr_{n+1}\Omega^{p+1}$. Let $\Phi^p_n$
be the kernel of this map. We have the exact sequence
\bea
&&
0\lar \Phi^p_n\lar \gr_n\,\Omega^p \stackrel{d}{\lar}
 d\gr_{n}\,\Omega^p\lar 0.
\non
\ena
We define the (graded version of) residue sheaf $R^p_n$ \cite{AH} (see also
appendix to chapter VII by Mumford in \cite{Z}) by
\bea
&&
R^p_n=\Phi^p_n/d\gr_{n-1}\,\Omega^{p-1},
\quad
n\geq 2.
\non
\ena
Notice that the support of $R^p_n$ is contained in $\{0\}$, since closed forms
are exact at a non-singular point of $\Theta$. 
In other words the de Rham complex may not be exact at a singular point of 
$\Theta$. In order to study the ${\cal D}$-module structure of abelian functions at the level of sheaves it is necessary to study $R^p_n$.
For $p=1$ we have

\begin{lemma}\label{lem-1}
(i) $\Phi^1_n= d\gr_{n-1}{\cal O}$ for $n\geq 2$.
\vskip2mm
\noindent
(ii) $\Phi^1_1\simeq \gr_0\,{\cal O}$.
\vskip2mm
\noindent
(iii) $d\gr_n\,{\cal O}\simeq \gr_n\,{\cal O}$ for $n\geq 1$.
\end{lemma}

To prove the lemma we analyze the defining equation of $\Theta$
near the singular point. The following proposition is well known
from the general theory of singularities \cite{A1,A2}. 
For the sake to be self-contained we give an elementary 
computational proof.

\begin{prop}\label{coord} There exists a local coordinate 
system $(z_1,z_2,z_3)$ near $0$ such that
\bea
&&
\sigma(u)=z_1^2+z_2^2+z_3^2.
\non
\ena
\end{prop}
\vskip2mm
\noindent
{\it Proof.}
Due to (ii) of Definition \ref{def-1} $\sigma(u)$ has the expansion of the form
\bea
&&
\sigma(u)=S(u)+\sum_{d<-6}f_d(u),
\quad
S(u)=u_1u_3-u_2^2-\frac{1}{3}u_1^3u_2+\frac{1}{45}u_1^6.
\label{g3-sigma}
\ena
Therefore $\sigma(u)$ can be written as
\bea
&&
\sigma(u)=u_1u_3-u_2^2+au_3^2+\sum_{d\geq 3} F_d(u),
\label{g3-exp}
\ena
where $F_d(u)$ is a homogeneous polynomial of degree $d$ with respect to
$\deg\, u_i=1$ $i=1,2,3$ and $a$ is some constant. We define $x_1$, $x_2$, $x_3$ by
\bea
&&
u_1+au_3=x_1+ix_2,
\quad
u_2=ix_3,
\quad
u_3=x_1-ix_2.
\non
\ena
Then
\bea
&&
\sigma(u)=x_1^2+x_2^2+x_3^2
+\left(\hbox{$\deg \geq 3$ terms in $x_1,x_2,x_3$}\right),
\non
\ena
where $\deg\, x_i=1$ $i=1,2,3$.
Therefore one can find holomorphic functions $G_1$, $G_2$, $G_3$ near 
$0$ and a constant $c$ such that
\bea
&&
\sigma(u)=x_1^2(1+G_1)+x_2^2(1+G_2)+x_3^2(1+G_3)+cx_1x_2x_3,
\non
\ena
and $G_i(0)=0$ for all $i$. Define the local coordinate $(X_1,X_2,X_3)$ by
\bea
&&
X_i=x_1(1+G_i)^{1/2},
\quad
i=1,2,3.
\non
\ena
Then
\bea
&&
x_1x_2x_3=GX_1X_2X_3,
\non
\ena
where $G=\prod_{i=1}^3(1+G_i)^{-1/2}$ can be considered as a 
holomorphic function of $X=(X_1,X_2,X_3)$ near $0$ such that $G(0)=1$.
With respect to variables $X$ we have
\bea
\sigma(u)&=&X_1^2+X_2^2+X_3^2+cGX_1X_2X_3,
\non
\\
&=&
(X_1+\frac{cG}{2}X_2X_3)^2+X_2^2(1-\frac{1}{4}c^2G^2X_3^2)+X_3^2.
\non
\ena
Take the local coordinates $(z_1,z_2,z_3)$ as
\bea
&&
z_1=X_1+\frac{cG}{2}X_2X_3,
\quad
z_2=X_2(1-\frac{1}{4}c^2G^2X_3^2)^{1/2},
\quad
z_3=X_3,
\non
\ena
we get
\bea
&&
\sigma(u)=z_1^2+z_2^2+z_3^2,
\non
\ena
which completes the proof of the proposition.
\qed
\vskip10mm

\noindent
{\it Proof of Lemma \ref{lem-1}}
\vskip2mm
\noindent
(i) The assertion is easily proved at a non-singular point of $\Theta$.
Let us prove the assertion at the singular point $0$. 
Take a local coordinate $(z_1,z_2,z_3)$ as in Proposition \ref{coord}.
Let 
\bea
&&
s=\sum_{i=1}^3\frac{a_i(z_1,z_2,z_3)}{\sigma^n}dz_i
\non
\ena
be a local section of $\gr_n\,\Omega^1$ around $0$. Notice that
\bea
&&
\frac{z_3^2h}{\sigma^n}=\frac{-(z_1^2+z_2^2)h}{\sigma^n}
\quad
\hbox{in $\gr_n\,\Omega^1$}
\non
\ena
for any holomorphic one form $h$ around $0$. Thus one can assume that
$a_i$ is linear in $z_3$:
\bea
&&
a_i=a_{i0}(z_1,z_2)+a_{i1}(z_1,z_2)z_3.
\non
\ena
Then the condition $ds=0$ in $\gr_{n+1}\,\Omega^2$ is equivalent to 
the following equations:
\bea
-z_2a_{10}+z_1a_{20}+z_3(-z_2a_{11}+z_1a_{21})&=&0,
\non
\\
z_1a_{30}+(z_1^2+z_2^2)a_{11}+z_3(-a_{10}+z_1a_{31})&=&0,
\non
\\
z_2a_{30}+(z_1^2+z_2^2)a_{21}+z_3(-a_{20}+z_2a_{31})&=&0.
\non
\ena
Then we have
\bea
&&
a_{10}=z_1a_{31},
\quad
a_{20}=z_2a_{31},
\quad
a_{11}=z_1b,
\quad
a_{21}=z_2b,
\quad
a_{30}=-(z_1^2+z_2^2)b,
\non
\ena
for some holomorphic function $b$ of $(z_1,z_2)$.
Consequently
\bea
&&
a_1=z_1(a_{31}+bz_3),
\quad
a_2=z_2(a_{31}+bz_3),
\quad
a_3=-(z_1^2+z_2^2)b+a_{31}z_3.
\non
\ena
Notice that
\bea
&&
a_3=z_3(a_{31}+bz_3)-b\sigma.
\non
\ena
Thus
\bea
&&
s=(a_{31}+bz_3)\frac{\sum_{i=1}^3z_idz_i}{\sigma^n}
=d\left(-\frac{1}{2(n-1)}\frac{a_{31}+bz_3}{\sigma^{n-1}}\right)
\quad
\hbox{in $\gr_{n}\,\Omega^1$},
\label{phi-1n}
\ena
since $n\geq 2$.

\vskip2mm
\noindent
(ii) By $(\ref{phi-1n})$ we have
\bea
\Phi^1_1={\cal O}\frac{d\sigma}{\sigma},
\non
\ena
around $0$. Then the map
\bea
gr_0\,{\cal O}&\lar& \Phi^1_1
\non
\\
F&\mapsto& Fd\log\,\sigma,
\non
\ena
gives an isomorphism.
\vskip2mm
\noindent
(iii) The proof is easy and we leave it to the reader.

\qed

Next we study the residue sheaf $R^p_n$ for $p\geq 2$.

\begin{lemma} Take a local coordinate system $(z_1,z_2,z_3)$ as in 
Proposition \ref{coord}. 
Then stalks at $0$ of $R^3_n$ and $R^2_n$ are described as
\vskip2mm
\noindent
(i) $\quad\displaystyle{(R^3_n)_0={\mathbf C} \varphi^3_n}$, 
$\quad \displaystyle{\varphi^3_n=\frac{dz_1\wedge dz_2\wedge dz_3}{\sigma^n}}$,
\vskip2mm
\noindent
(ii) $\quad\displaystyle{(R^2_n)_0={\mathbf C} \varphi^2_n}$, 
$\quad \displaystyle{\varphi^2_n=\frac{z_1dz_2\wedge dz_3+z_2dz_3\wedge dz_1+
z_3dz_1\wedge dz_2}{\sigma^n}}$.
\end{lemma}
\vskip2mm
\noindent
{\it Proof.} (i) The assertion follows from
\bea
&&
d\left(
\frac{a_1dz_2\wedge dz_3+a_2dz_3\wedge dz_1+
a_3dz_1\wedge dz_2}{\sigma^n}
\right)
=
-2n\frac{\sum_{i=1}^3z_ia_i}{\sigma^{n+1}}dz_1\wedge dz_2\wedge dz_3,
\non
\ena
with $a_i$ being holomorphic a function at $0$.

\vskip2mm
\noindent
(ii) Let
\bea
s&=&\frac{a_1dz_2\wedge dz_3+a_2dz_3\wedge dz_1+
a_3dz_1\wedge dz_2}{\sigma^n},
\non
\\
a_i&=&a_{i0}(z_1,z_2)+a_{i1}(z_1,z_2)z_3,
\quad
i=1,2,
\non
\ena
be a local section of $\gr_n\,\Omega^2$ at $0$. Then $ds=0$ in 
$\gr_{n+1}\,\Omega^3$ is equivalent to
\bea
a_{10}z_1+a_{20}z_2-a_{31}(z_1^2+z_2^2)&=&0,
\non
\\
a_{11}z_1+a_{21}z_2+a_{30}&=&0.
\non
\ena
It follows that there exists a holomorphic function $b(z_1,z_2)$ at $0$ 
such that
\bea
a_{10}-a_{31}z_1&=&bz_2,
\non
\\
a_{20}-a_{31}z_2&=&-bz_1,
\non
\ena
and consequently
\bea
a_1&=&a_{31}z_1+bz_2+a_{11}z_3,
\non
\\
a_{2}&=&a_{31}z_2-bz_1+a_{21}z_3,
\non
\\
a_3&=&-a_{11}z_1-a_{21}z_2+a_{31}z_3.
\non
\ena
Then we have
\bea
s&=&a_{31}\varphi^2_n+
\frac{1}{2(1-n)}
d\left(
\frac{a_{21}dz_1-a_{11}dz_2+bdz_3}{\sigma^{n-1}}
\right).
\non
\ena
Therefore $\Phi^2_n$ modulo $d(\gr_{n-1}\,\Omega^1)$ is represented
by forms of the form $a\varphi^2_n$ with $a$ being holomorphic at $0$.
Notice that
\bea
&&
z_ia\varphi^2_n=
\frac{1}{2(1-n)}
d\left(
a\frac{z_{i+2}dz_{i+1}-z_{i+1}dz_{i+2}}{\sigma^{n-1}}
\right),
\non
\ena
where the index of $z_i$ should be read modulo $3$. Thus $R^2_n$ is
represented by elements of ${\mathbf C}\varphi^2_n$.
Using
\bea
&&
d\left(
\frac{\sum_{i=1}^3a_idz_i}{\sigma^n}
\right)
\non
\\
&&
=-2n
\frac{
(z_2a_3-z_3a_2)dz_2\wedge dz_3
+(z_3a_1-z_1a_3)dz_3\wedge dz_1
+(z_1a_2-z_2a_1)dz_1\wedge dz_2
}
{\sigma^{n+1}},
\non
\ena
one can easily check that $\varphi^2_n$ is not zero in $R^2_n$.
\qed
\vskip5mm

Since $d\Phi^2_n=0$ in $\gr_{n+1}\,\Omega^3$, the map
\bea
&&
d:R^2_n\lar \gr_{n}\,\Omega^3/d\gr_{n-1}\,\Omega^2=R^3_n
\label{r2-r3}
\ena
is well defined.

\begin{lemma}\label{lem-r2-r3} The map (\ref{r2-r3}) is an isomorphism of
${\cal O}$-modules.
\end{lemma}
\vskip2mm
\noindent
{\it Proof.} The lemma follows from
\bea
&&
d\varphi^2_n=3\varphi^3_n.
\non
\ena
\qed

\section{Finite Generation of $A$ over ${\cal D}$}
In this section we study the differential structure of the cohomology
groups $H^0(J(X),\gr_n\,\Omega^p)$ and prove that $A$ is a finitely generated
${\cal D}$-module.

We use the following vanishing theorem.

\begin{theorem}\label{vanish}
\par
\vskip2mm
\noindent
(i) $\displaystyle{H^i(J(X),{\cal O}(n))=0}$ for $n\geq 1$, $i\geq 1$.
\vskip2mm
\noindent
(ii) $\displaystyle{H^i(J(X),\gr_n\,{\cal O})=0}$ for $n\geq 2$, $i\geq 1$.
\vskip2mm
\noindent
(iii) $\displaystyle{H^i(J(X),\gr_1\,{\cal O})\simeq H^{i+1}(J(X),{\cal O})}$ 
for $i\geq 0$.
\vskip2mm
\noindent
(iv) $\displaystyle{H^i(J(X),\gr_0\,{\cal O})\simeq\left\{
\begin{array}{rl}
H^i(J(X),{\cal O})&\quad i\leq 2,\\
0&\quad i>2.\\
\end{array}
\right.
}
$
\end{theorem}

The assertion (i) is due to Mumford \cite{M0} and (ii), (iii) follows from it
using the exact sequence
\bea
&&
0\lar {\cal O}(n-1)\lar {\cal O}(n) \lar \gr_n\, {\cal O} \lar 0.
\label{seq}
\ena

Notice that
\bea
&&
H^i(J(X),\gr_n\,\Omega^p)\simeq H^i(J(X),\gr_n\,{\cal O})\otimes H^0(J(X),\Omega^p),
\non
\ena
since $\Omega^p$ is a free ${\cal O}$ module.

By the definition
\bea
&&
A(n)=H^0(J(X),{\cal O}(n)).
\non
\ena
Due to (i) of Theorem \ref{vanish} and (\ref{seq}) we have
\bea
&&
H^0(J(X),\gr_n\,{\cal O})=\gr^P_n\,A
\quad
\hbox{for $n\geq 2$}.
\non
\ena

We shall study the ${\cal D}$-module structure of $\gr^P\,A$.
In the generic case where $\Theta$ is non-singular $\gr^P\,A$ is
finitely generated over ${\cal D}$. It means that functions with
n-th order poles can be obtained by differentiating functions with
(n-1)-st order poles for all sufficiently large $n$ \cite{CN}.
The following proposition shows that it does not hold in the
 present case. This is due to the existence of the singularity
of $\Theta$ (recall that $R^p_n$ vanishes if $\Theta$ is non-singular).
More precisely (i) of the proposition implies that there exists,
 up to constant multiples, a "missing function" in $A(n)$ such that
a linear combination of derivatives of them again belongs to $A(n)$.
The statement (ii) of the proposition shows that there are functions
of $A(n)$ such that a linear combination of derivatives of them again 
belongs to $A(n)$. In Proposition \ref{gen-4} we prove that the missing
function in $A(n)$ is obtained from this linear combination in $A(n)$.

\begin{prop}\label{fg}
\par
\vskip5mm
\noindent
(i) $\displaystyle{
\frac{H^0(J(X),\gr_n\,\Omega^3)}{dH^0(J(X),\gr_{n-1}\Omega^2)}
\simeq H^0(J(X),R^3_n)
\,}$ for $n\geq 5$.
\vskip2mm
\noindent
(ii) $\displaystyle{\frac{\ker\left(d:H^0(J(X),\gr_n\,\Omega^2)\lar 
H^0(J(X),\gr_{n+1}\,\Omega^3)\right)}{dH^0(J(X),\gr_{n-1}\,\Omega^1)}
\simeq H^0(J(X),R^2_n)\,}$ for $n\geq 4$.
\end{prop}
\vskip5mm
\noindent
{\it Proof.} 
The cohomology sequence of
\bea
&&
0\lar d\gr_{n-1}\,\Omega^2 \lar \gr_n\,\Omega^3 \lar R^3_n \lar 0,
\quad
n\geq 2
\label{short-1}
\ena
gives 
\bea
&&
H^i(J(X),d\gr_{n-1}\,\Omega^2)=0,
\quad
n\geq 2,\quad i\geq 2,
\label{iso-1}
\ena
and the exact sequence
\bea
&&
0 \lar H^0(J(X),d\gr_{n-1}\,\Omega^2) \lar H^0(J(X),\gr_n\,\Omega^3) \lar
H^0(J(X),R^3_n)
\non
\\
&&
\lar 
H^1(J(X),d\gr_{n-1}\,\Omega^2) \lar 0,
\label{ex-1}
\ena
by Theorem \ref{vanish} (i) and $H^i(J(X),R^3_n)=0$, $i\geq 1$.
Let us first prove
\bea
&&
H^1(J(X),d\gr_{n-1}\,\Omega^2)=0
\quad
\hbox{for $n\geq 4$}.
\label{aim-1}
\ena
The exact sequence
\bea
&&
0 \lar \Phi^2_{n-1} \lar \gr_{n-1}\,\Omega^2 
\stackrel{d}{\lar} d\gr_{n-1}\,\Omega^2 \lar 0.
\label{short-2}
\ena
implies the isomorphism
\bea
&&
H^i(J(X),d\gr_{n-1}\,\Omega^2)\simeq H^{i+1}(J(X),\Phi^2_{n-1}),
\quad
i\geq 1,\quad n\geq 3,
\label{iso-2}
\ena
and the exact sequence
\bea
&&
0 \lar H^0(J(X),\Phi^2_{n-1}) \lar H^0(J(X),\gr_{n-1}\,\Omega^2) 
\lar H^0(J(X),d\gr_{n-1}\,\Omega^2) 
\non
\\
&&
\lar H^1(J(X),\Phi^2_{n-1}) \lar 0,
\quad
n\geq 3.
\label{ex-2}
\ena
Similarly, by the exact sequence,
\bea
&&
0 \lar d\gr_{n-2}\,\Omega^1 \lar \Phi^2_{n-1} \lar R^2_{n-1} \lar 0,
\label{short-3}
\ena
we get
\bea
&&
H^i(J(X),d\gr_{n-2}\,\Omega^1)\simeq H^i(J(X),\Phi^2_{n-1}),
\quad
i\geq 2,\quad n\geq 3,
\label{iso-3}
\ena
and the exact sequence
\bea
&&
0 \lar H^0(J(X),d\gr_{n-2}\,\Omega^1) \lar H^0(J(X),\Phi^2_{n-1}) \lar H^0(J(X), R^2_{n-1})
\non
\\
&&
\lar H^1(J(X),d\gr_{n-2}\,\Omega^1) \lar H^1(J(X),\Phi^2_{n-1}) \lar 0,
\quad
n\geq 3.
\label{ex-3}
\ena
Considering
\bea
&&
0 \lar \Phi^1_{n-2} \lar \gr_{n-2}\,\Omega^1 
\stackrel{d}{\lar} d\gr_{n-2}\,\Omega^1 \lar 0,
\label{short-4}
\ena
we have
\bea
&&
H^i(J(X),d\gr_{n-2}\, \Omega^1)\simeq H^{i+1}(J(X),\Phi^1_{n-2}),
\quad
i\geq 1,\quad n\geq 4,
\label{iso-4}
\ena
and the exact sequence
\bea
&&
0 \lar H^0(J(X),\Phi^1_{n-2}) \lar H^0(J(X),\gr_{n-2}\, \Omega^1) 
\stackrel{d}{\lar} 
H^0(J(X),d\gr_{n-2}\, \Omega^1) 
\non
\\
&&
\lar H^1(J(X), \Phi^1_{n-2}) \lar 0,
\quad
n\geq 4.
\label{ex-4}
\ena
By Lemma \ref{lem-1} (i), (iii) 
\bea
&&
\Phi^1_{n-2}\simeq \gr_{n-3}\,{\cal O},
\quad
n\geq 4.
\label{iso-5}
\ena
Using (\ref{iso-2}), (\ref{iso-3}), (\ref{iso-4}), (\ref{iso-5}) we get, 
for $n\geq 4$,
\bea
&&
H^1(J(X),d\gr_{n-1}\,\Omega^2)
\simeq 
H^2(J(X),\Phi^2_{n-1})
\simeq 
H^2(J(X),d\gr_{n-2}\,\Omega^1)
\non
\\
&&
\simeq
H^3(J(X),\Phi^1_{n-2})
\simeq
H^3(J(X),\gr_{n-3}\,{\cal O}).
\label{iso-5+}
\ena
It vanishes, since the support of $\gr_{n-3}\,{\cal O}$ 
is contained in $\Theta$ and $\dim\,\Theta=2$.

Next we prove
\bea
&&
H^0(J(X),d\gr_{n-1}\,\Omega^2)=dH^0(J(X),\gr_{n-1}\,\Omega^2).
\label{aim-2}
\ena
We have, by (\ref{iso-4}), (\ref{iso-5}) and Theorem \ref{vanish} (ii),
\bea
&&
H^1(J(X),d\gr_{n-2}\,\Omega^1)\simeq H^2(J(X),\Phi^1_{n-2})\simeq 
H^2(J(X),\gr_{n-3}\,{\cal O})=0
\quad
n\geq 5,
\label{iso-6}
\ena
and, by (\ref{ex-3}), 
\bea
&&
H^1(J(X),\Phi^2_{n-1})=0,
\quad
n\geq 5.
\label{iso-7}
\ena
Then the equation (\ref{aim-2}) follows from (\ref{ex-2}) and 
 claim (i) follows from (\ref{ex-1}).

\vskip2mm
\noindent
(ii) Notice that
\bea
&&
\ker\left(
d: H^0(J(X),\gr_n\,\Omega^2) \lar H^0(J(X),\gr_{n+1}\,\Omega^3)
\right)
=H^0(J(X),\Phi^2_n).
\label{iso-7+}
\ena
We have
\bea
&&
\frac{H^0(J(X),\Phi^2_n)}{H^0(J(X),d\gr_{n-1}\,\Omega^1)}
\simeq 
H^0(J(X),R^2_n),
\quad
n\geq 4,
\label{iso-8}
\ena
by (\ref{ex-3}), (\ref{iso-6}) and
\bea
&&
\frac{H^0(J(X),d\gr_{n-1}\,\Omega^1)}{dH^0(J(X),\gr_{n-1}\,\Omega^1)}
\simeq 
H^1(J(X),\Phi^1_{n-1}),
\quad
n \geq 3,
\label{iso-9}
\ena
by (\ref{ex-4}).
Then the assertion (ii) 
follows from (\ref{iso-8}) and (\ref{iso-9}) 
using
\bea
&&
H^1(J(X),\Phi^1_{n-1})\simeq H^1(J(X),\gr_{n-2}\,{\cal O})=0,
\quad
n\geq 4.
\non
\ena
\qed

\begin{prop}\label{gen-4}
The ${\cal D}$ module $A$ is generated by $A(4)$.
\end{prop}
\vskip2mm
\noindent
{\it Proof.} 
By Proposition \ref{fg} we have, for $n\geq 5$,
\bea
A(n)&=&V^n_1+{\mathbb C} F_n+A(n-1),
\non
\\
V^n_1&=&\sum_{i=1}^3\partial_i A(n-1),
\non
\ena
where $F_n$ is an element of $A(n)$ such that
\bea
&&
\varphi^3_n=F_ndu_1\wedge du_2\wedge du_3
\quad
\hbox{in $(R^3_n)_0$}.
\non
\ena
It means, in particular, that
\bea
&&
H^0(J(X),R^3_n)={\mathbb C}F_n.
\non
\ena

\begin{lemma}\label{lem-41}
We have, for $n\geq 5$,
\bea
&&
F_n\in \sum_{i=1}^3\partial_i V^n_1+V^n_1+A(n-1).
\non
\ena
\end{lemma}
\vskip2mm
\noindent
{\it Proof.}
By Lemma \ref{lem-r2-r3}
\bea
&&
dH^0(J(X),R^2_n)=H^0(J(X),R^3_n).
\non
\ena
Claim (ii) of Proposition \ref{fg} implies that an element of 
$H^0(J(X),R^2_n)\simeq {\mathbb C}$ is represented by an element 
of $H^0(J(X),\Omega^2(n))$. Let us take elements $f_i$, $i=1,2,3$ 
of $A(n)$ such that the two form
\bea
&&
f_1du_2\wedge du_3+f_2du_3\wedge du_1+f_3du_1\wedge du_2
\non
\ena
is a basis of $H^0(J(X),R^2_n)$ and it coincides with 
$\varphi^2_n$ in $(R^2_n)_0$.
Then $F_n$ can be written in a form
\bea
F_n=\sum_{i=1}^3\frac{\partial f_i}{\partial u_i}+
\sum_{i=1}^3\frac{\partial \tilde{f}_i}{\partial u_i}
+G_{n-1},
\non
\ena
for some $\tilde{f}_i,G_{n-1}\in A(n-1)$.
Since $f_i\in A(n)$ one can write
\bea
&&
f_i=g_i+c_iF_n+h_i,
\label{ci=0}
\ena
where $c_i$ is a constant and
\bea
&&
g_i\in V^n_1,
\quad
h_i\in A(n-1).
\non
\ena
\vskip5mm

\begin{lemma}\label{lem-42}
$\displaystyle{c_i=0}$.
\end{lemma}
\vskip2mm
\noindent
{\it Proof.}
Multiply $\sigma^n$ to Equation (\ref{ci=0})
and set $u=(0,0,0)$. By changing the coordinate to $(z_1,z_2,z_3)$ as 
in Proposition \ref{coord} around $0$, we see that the right hand side
becomes $c_i$. On the other hand the left hand side is zero. Because
$\varphi^2_n\sigma^n$ and elements of $V^n_1\sigma^n$ vanish at $0$.
Thus $c_i=0$.
\qed
\vskip5mm

By this lemma we have
\bea
&&
f_i\in V^n_1+A(n-1).
\non
\ena
Lemma \ref{lem-41} follows from this.
\qed
\vskip5mm

By Lemma \ref{fg} we have
\bea
&&
A(n)\subset {\cal D}A(n-1)
\quad
\hbox{for $n\geq 5$}.
\non
\ena
Thus $A$ is generated by $A(4)$ over ${\cal D}$.
\qed
\vskip5mm

Finally notice the following corollary of Proposition \ref{fg}.

\begin{cor}\label{not-fg} As a ${\cal D}$-module 
$\gr^P\,A$ is not finitely generated.
\end{cor}

On the other hand $\gr^{KP}\,A$ is finitely generated over 
${\cal D}$ as we shall see later (Theorem \ref{main-th}).
\vskip2mm

\noindent
\begin{remark} For $g=2$ $\gr^P\,A$ is finitely generated
as proved in \cite{CN}.
\end{remark}

\section{Cohomology of Affine Jacobian}
In this section we briefly recall the results on a description of 
the cohomology group $H^3(J(X)\backslash\Theta,{\mathbb C})$ \cite{N1}.

By the algebraic de Rham theorem we have the isomorphism
\bea
&&
H^3(J(X)\backslash\Theta,{\mathbb C})\simeq 
\left(A/\sum_{i=1}^3 \partial_i A\right)\otimes \wedge^3 T^\ast.
\non
\ena

\begin{theorem}\label{sing-cohom}\cite{N1}
There is an isomorphism
\bea
&&
H^3(J(X)\backslash\Theta,{\mathbb C})\simeq W^3,
\label{h3-w3}
\ena
where $W^3$ is given by (\ref{wk}) with $g=3$ and $k=3$.
The composition of maps $\ev$ and the projection 
$A\rightarrow A/\sum_{i=1}^3 \partial_i A$ gives the isomorphism
\bea
&&
W^3\lar \left(A/\sum_{i=1}^3 \partial_i A\right)\otimes \wedge^3 T^\ast.
\non
\ena
\end{theorem}

It follows from Theorem \ref{sing-cohom} that
\bea
&&
\dim\, H^3(J(X)\backslash\Theta,{\mathbb C})=14,
\non
\ena
and that $A/\sum_{i=1}^3 \partial_i A$ is generated over ${\mathbb C}$
by
\bea
&&
1,
\quad
(i;j)=\zeta_{ij},
\quad
(ij;kl),
\quad
(123;123).
\label{c-gen}
\ena

\begin{prop}\label{cohomology-basis}
A basis of $A/\sum_{i=1}^3 \partial_i A$ is given by
\bea
&&
1,
\quad
\zeta_{ij} 
\quad
(1\leq i\leq j\leq 3),
\non
\\
&&
(12;12),
\quad (12;13),
\quad (12;23),
\quad (13;13),
\quad (13;23),
\quad (23;23),
\non
\\
&&
(123;123).
\label{c-basis}
\ena
\end{prop}
\vskip2mm
{\it Proof.} Notice that, by the definition,
 $(i_1,...,i_k;j_1,...,j_k)$ is skew symmetric
in $i_1,...,i_k$ and $j_1,...,j_k$ respectively and satisfies the symmetry 
relation
\bea
&&
(i_1,...,i_k;j_1,...,j_k)=(j_1,...,j_k;i_1,...,i_k).
\non
\ena
It follows that any element of (\ref{c-gen}) is a constant multiple 
of an element in (\ref{c-basis}). Since the number of elements in 
(\ref{c-basis}) is $14$, they are linearly independent.
\qed

\begin{lemma}\label{ord-zeta}
(i) $\ord\,\zeta_{ij}=2$.
\vskip2mm
\noindent
(ii) $\ord\,(ij;kl)\leq 3$.
\vskip2mm
\noindent
(iii) $\ord\,(ijk;lmn)\leq 4$.
\end{lemma}
\vskip2mm
\noindent
{\it Proof.} It is proved in Lemma 8.3 of \cite{CN} that
\bea
&&
\ord\,(i_1,...,i_k;j_1,...,j_k)\leq k+1.
\label{pole-ord}
\ena
The assertions (ii), (iii) follow from this. Let us prove (i).
Obviously $\ord\,\zeta_{ij}\leq 2$. If
$\ord\,\zeta_{ij}<2$, then $\zeta_{ij}$ is a constant.
Because $A(1)={\mathbb C}$. It contradicts Proposition 
\ref{cohomology-basis} which claims, in particular,
 the linear independence of $\{1,\zeta_{ij}\}$.  
\qed

\section{Baker's Addition Formula}
In order to describe a basis of $A(2)$ in terms of the basis
of the cohomology group of the affine Jacobian given 
in Proposition \ref{cohomology-basis} 
we use the addition formula of the sigma function \cite{B1,B2,BEL1}.

Let $u=(u_1,u_2,u_3)$ and $v=(v_1,v_2,v_3)$.
The addition formula for the $g=3$ hyperelliptic sigma function due to
Baker \cite{B2} is 
\bea
&&
\frac{\sigma(u+v)\sigma(u-v)}{\sigma(u)^2\sigma(v)^2}
=
(\wp_{13}(v)-\wp_{13}(u))(\wp_{22}(v)-\wp_{22}(u))-(\wp_{13}(v)-\wp_{13}(u))^2
\non
\\
&&
\qquad\qquad\qquad\qquad\quad
-(\wp_{23}(v)-\wp_{23}(u))(\wp_{12}(v)-\wp_{12}(u))
\non
\\
&&
\qquad\qquad\qquad\qquad\quad
+
(\wp_{33}(v)-\wp_{33}(u))(\wp_{11}(v)-\wp_{11}(u)).
\label{add-1}
\ena
We set
\bea
&&
\Sym(F(u)G(v))=F(u)G(v)+F(v)G(u).
\non
\ena
Then (\ref{add-1}) is rewritten as
\bea
&&
\frac{\sigma(u+v)\sigma(u-v)}{\sigma(u)^2\sigma(v)^2}
=
\Sym\left((13;13)(u)-(12;23)(u))\cdot1\right)
-\Sym\left(\wp_{13}(u)\wp_{22}(v)\right)
\non
\\
&&
\qquad\qquad\qquad\qquad\quad
+\Sym\left(\wp_{13}(u)\wp_{13}(v)\right)
+\Sym\left(\wp_{12}(u)\wp_{23}(v)\right)
\non
\\
&&
\qquad\qquad\qquad\qquad\quad
-\Sym\left(\wp_{11}(u)\wp_{33}(v)\right).
\label{add-2}
\ena

\begin{cor}\label{add-cor1}
$\displaystyle{\ord\,\left((13;13)-(12;23)\right)=2}$.
\end{cor}
\vskip2mm
\noindent
{\it Proof.} Consider the equation (\ref{add-2}) as that for functions
of $u$. Then all terms other than $(13;13)(u)-(12;23)(u)$ have poles of 
order less than or equal to two. Thus the order of poles of 
$(13;13)(u)-(12;23)(u)$
is at most two. If the order of poles is less than two 
$(13;13)(u)-(12;23)(u)$ becomes a constant. 
It contradicts the linear independence of 
$1$, $(13;13)$, $(12;23)$ given by Proposition 
\ref{cohomology-basis}.\qed

\begin{cor}\label{basis-a2}
(i)
$\displaystyle{
A(2)={\mathbb C}1\oplus \oplus_{1\leq i\leq j\leq 3}{\mathbb C}\wp_{ij}
\oplus{\mathbb C}\left((13;13)-(12;23)\right).
}$
\vskip5mm
\noindent
(ii) $\displaystyle{\gr^P_2\, A=\oplus_{1\leq i\leq j\leq 3}{\mathbb C}\wp_{ij}}
\oplus{\mathbb C}\left((13;13)-(12;23)\right).$
\end{cor}
\vskip2mm
\noindent
{\it Proof.} By Corollary \ref{add-cor1} $(13;13)-(12;23)\in A(2)$.
By Proposition \ref{cohomology-basis} elements appeared in (i)
are linearly independent. Since $\dim\, A(2)=8$, they forms a basis of $A(2)$.
The assertion (ii) follows from (i).
\qed

\section{Abelian Functions of Order Three}
In this section we study $A(3)$.

\begin{prop}\label{order=3}
(i) $\displaystyle{
\frac{H^0(J(X),\gr_3\,\Omega^3)}{H^0(J(X),d\gr_2\,\Omega^2)}\simeq H^0(J(X),R^3_3).
}$
\vskip3mm
\noindent
(ii) $\displaystyle{
\frac{H^0(J(X),d\gr_2\,\Omega^2)}{dH^0(J(X),\gr_2\,\Omega^2)}\simeq H^1(J(X),\Phi^2_2).
}$
\vskip3mm
\noindent
(iii) $\displaystyle{\dim\, H^1(J(X),\Phi^2_2)=5}$.
\end{prop}
\vskip2mm
\noindent
{\it Proof.} (i) By (\ref{ex-1}) and (\ref{iso-2}) with $n=3$ it is sufficient
to prove
\bea
&&
H^2(J(X),\Phi^2_2)=0.
\label{71}
\ena
Since $\Phi^1_1\simeq \gr_0\,{\cal O}$ by Lemma \ref{lem-1} (ii), $H^3(J(X),\Phi^1_1)=0$ by Theorem \ref{vanish} (iv). The long cohomology exact sequence
of (\ref{short-4}) with $n=3$ is
\bea
&&
0\lar H^0(J(X),\Phi^1_1) \lar H^0(J(X),\gr_1\,\Omega^1) \lar H^0(J(X),d\gr_1\,\Omega^1) \lar
\cdots
\non
\\
&&
\lar H^2(J(X),\Phi^1_1) \lar H^2(J(X),\gr_1\,\Omega^1) \lar H^2(J(X),d\gr_1\,\Omega^1) \lar 0.
\non
\ena

\begin{lemma}\label{lem-71}
For $0\leq i\leq 2$ the map $\alpha: H^i(J(X),\Phi^1_1) \lar H^i(J(X),\gr_1\,\Omega^1)$
is injective.
\end{lemma}
\vskip2mm
\noindent
{\it Proof.} By (iii), (iv) of Theorem \ref{vanish} 
we have, for $0\leq i\leq 2$,
\bea
&&
H^i(J(X),\Phi^1_1)\simeq H^i(J(X),{\cal O})\simeq \wedge^i {\bar{T}}^\ast,
\non
\\
&&
H^i(J(X),\gr_1\,\Omega^1)\simeq H^{i+1}(J(X),{\cal O})\otimes H^0(J(X),\Omega^1)
\simeq \wedge^{i+1}{\bar{T}}^\ast\wedge T^\ast,
\non
\ena
where ${\bar{T}}^\ast=\sum_{i=0}^g {\mathbb C}d{\bar u}_i$ and
$d{\bar u}_i$ is the complex conjugate of $du_i$.
Then the lemma is proved in a similar manner to Lemma 4.6 in \cite{CN}
using the representation theory of $sl_2$.
\qed

By Lemma \ref{lem-71} we have the exact sequences;
\bea
&&
0 \lar H^i(J(X),\Phi^1_1) \lar H^i(J(X),\gr_1\,\Omega^1) \lar H^i(J(X),d\gr_1\,\Omega^1) 
\lar 0,
\label{72}
\ena
for $i \leq 2$.
Then
\bea
\dim\,H^2(J(X),d\gr_1\,\Omega^1)&=&\dim\,H^2(J(X),\gr_1\,\Omega^1)-\dim\,H^2(J(X),\Phi^1_1)
\non
\\
&=&
\dim\,H^3(J(X),\Omega^1)-\dim\, H^2(J(X),{\cal O})
\non
\\
&=&0.
\non
\ena
Here we use Theorem \ref{vanish} (iii), (iv). Thus $H^2(J(X),d\gr_1\,\Omega^1)=0$.
By (\ref{iso-3}) we have (\ref{71}).
\vskip2mm

\noindent
(ii) This follows from (\ref{ex-2}).
\vskip2mm

\noindent
(iii) Consider the exact sequence (\ref{ex-3}).
By (\ref{72}) 
\bea
&&
\dim\, H^1(J(X),d\gr_1\,\Omega^1)=6.
\non
\ena
Thus $\dim\, H^1(J(X),\Phi^2_2)=5$ or $6$, since $\dim\, H^0(J(X),R^2_2)=1$.
We prove that the latter is impossible.

Suppose that $\dim\, H^1(J(X),\Phi^2_2)=6$. Then
\bea
&&
\frac{H^0(J(X),\Phi^2_2)}{H^0(J(X),d\gr_1\,\Omega^1)}\simeq H^0(J(X),R^2_2).
\non
\ena
By (\ref{iso-7+}) and
\bea
&&
H^0(J(X),\gr_2\,\Omega^2)\simeq \frac{H^0(J(X),\Omega^2(2))}{H^0(J(X),\Omega^2(1))},
\non
\ena
there is an element $w$ of $H^0(J(X),\Omega^2(2))$ such that it is contained in 
$H^0(J(X),\Phi^2_2)$ and its image in $H^0(J(X),R^2_2)$ becomes a basis of $H^0(J(X),R^2_2)$. 
By Lemma \ref{lem-r2-r3} $dw$ becomes a basis of $H^0(J(X),R^3_2)$.
In particular
$dw$ is a non-zero element of $H^0(J(X),\Omega^3(2))$.
By Corollary \ref{basis-a2} a basis of $H^0(J(X),\Omega^3(2))$ is given by a subset 
of a basis of the cohomology group 
$H^3(J(X)\backslash\Theta,{\mathbb C})$.
Thus they are linearly independent modulo exact forms. Then $dw=0$ 
as an element of $H^0(J(X),\Omega^3(2))$. This is a contradiction.
Thus the assertion (iii) is proved.
\qed
\vskip5mm

Let us find a basis of the space appeared in Proposition \ref{order=3} (ii).

\begin{lemma}\label{lem-72}
For any $i,j,k,l$ $(ij;kl)\,\dum$ is in 
$H^0(J(X),d\gr_2\,\Omega^2)$.
\end{lemma}
\vskip2mm
\noindent
{\it Proof.}
Notice that
\bea
&&
d\zeta_i \wedge d\zeta_j \wedge du_k=
\frac{1}{2}d(\zeta_id\zeta_j\wedge du_k-\zeta_jd\zeta_i\wedge du_k).
\non
\ena
Here
\bea
&&
\ord\,(\zeta_id\zeta_j\wedge du_k-\zeta_jd\zeta_i\wedge du_k)\leq 2.
\non
\ena
In fact
\bea
\zeta_id\zeta_j\wedge du_k-\zeta_jd\zeta_i\wedge du_k
=\sum_{l\neq k}(\zeta_i\zeta_{jl}-\zeta_j\zeta_{il})du_l \wedge du_k,
\non
\ena
and
\bea
&&
\zeta_i\zeta_{jl}-\zeta_j\zeta_{il}=
\frac{\sigma_i\sigma_{jl}-\sigma_j\sigma_{il}}{\sigma^2},
\non
\ena
where
\bea
&&
\sigma_i=\partial_i\sigma, 
\qquad
\sigma_{ij}=\partial _i\partial_j\sigma.
\non
\ena
Thus the lemma is proved.
\qed
\vskip5mm

For the sake of simplicity we set
\bea
&&
v^0=(13:13)-(12:23)
\quad
v^1:=(12;12),
\quad
v^2:=(12;13),
\non
\\
&&
v^3:=(12;23),
\quad
v^4:=(13;23),
\quad
v^5:=(23;23).
\label{v-i}
\ena
Notice that
\bea
&&
(13;13)=(12;23)
\qquad
\hbox{in $\gr^P_3\,A$},
\non
\ena
by Corollary \ref{add-cor1}.

\begin{cor}\label{cor-71}
We have
\bea
&&
\frac{H^0(J(X),d\gr_2\,\Omega^2)}{dH^0(J(X),\gr_2\,\Omega^2)}
=\oplus_{i=1}^5 {\mathbb C}v^i\dum.
\non
\ena
\end{cor}
\vskip2mm
\noindent
{\it Proof.}
The left hand side is five dimensional by 
Proposition \ref{order=3} (ii), (iii).
Thus it is sufficient to prove the linear independence of $\{v_i\dum\}$ 
in the space of the left hand side.

Suppose that
$\sum c_iv_i\dum=0$ in this space. It means that there is an element
$w$ in $H^0(J(X),\Omega^2(2))$ such that
\bea
&&
\sum c_iv_i\dum -dw\in H^0(J(X),\Omega^3(2)).
\non
\ena
It implies $c_i=0$ for any $i$.
Since $\{\,v^i\dum\,|1\leq i\leq 5\,\}$ and the basis of 
$H^0(J(X),\Omega^3(2))=A(2)\dum$ given in (i) of Corollary \ref{basis-a2}
constitute a part of a basis of $H^3(J(X)\backslash\Theta,{\mathbb C})$.
\qed
\vskip5mm

We set
\bea
&&
w_{i_1...i_n}=\partial_{i_1}\cdots\partial_{i_n}w,
\non
\ena
for $w\in A$.

Let $F_3$ be an element of $A(3)$ such that $F_3\dum$ is a basis
of $H^0(R^3_3)$.

\begin{cor}\label{cor-72}
We have that
\bea
&&
\gr^P_3\,A=\oplus_{1\leq i\leq j\leq k\leq 3}{\mathbb C}\wp_{ijk}(u)
\oplus\oplus_{i=1}^3{\mathbb C}v^0_i
\oplus
\oplus_{i=1}^5{\mathbb C}v^i
\oplus {\mathbb C}F_3.
\non
\ena
\end{cor}
\vskip2mm
\noindent
{\it Proof.} By Proposition \ref{order=3} and Corollary \ref{cor-71},
$\gr^P_3\,A$ is generated by
\bea
&&
\wp_{ijk}(u) (1\leq i\leq j\leq k\leq 3),
\quad
v^0_i (1\leq i\leq 3),
\quad
v^i (1\leq i\leq 5),
\quad
F_3.
\non
\ena
The number of those elements is $19$. While 
$\dim\,\gr^P_3\,A=3^3-2^3=19$. Thus the above set of 
elements is linearly independent.
\qed

\section{Abelian Functions of Order Four}
In this section we study the space $A(4)$ and determine a minimal set of 
generators of the ${\cal D}$-module $A$.

\begin{prop}\label{prop-81}
We have the isomorphism
$$
\displaystyle{\frac{H^0(J(X),\gr_4\,\Omega^3)}{dH^0(J(X),\gr_3\,\Omega^2)}
\simeq H^0(J(X),R^3_4)}.
$$
\noindent
\end{prop}
\vskip2mm
\noindent
{\it Proof.} As proved in (\ref{ex-1}), (\ref{aim-1}), (\ref{iso-5+})
\bea
&&
\frac{H^0(J(X),\gr_4\,\Omega^3)}{H^0(J(X),d\gr_3\,\Omega^2)}
\simeq H^0(J(X),R^3_4).
\non
\ena
Let us prove
\bea
&&
H^0(J(X),d\gr_3\,\Omega^2)=dH^0(J(X),\gr_3\,\Omega^2).
\label{8-1}
\ena
By (\ref{ex-2}) it is equivalent to
\bea
&&
H^1(J(X),\Phi^2_3)=0.
\label{8-2}
\ena
In the exact sequence (\ref{ex-3}) with $n=4$
\bea
&&
H^1(J(X),d\gr_2\,\Omega^1)\simeq H^2(J(X),\Phi^1_2) 
\simeq H^2(J(X),\gr_1\,{\cal O}) \simeq H^3(J(X),{\cal O})\simeq {\mathbb C},
\non
\\
&&
H^0(J(X),R^2_3)\simeq {\mathbb C},
\non
\ena
by (\ref{iso-4}), Lemma \ref{lem-1} (i), (iii) and Theorem \ref{vanish} (iii).
Thus $\dim\, H^1(J(X),\Phi^2_3)=0$ or $1$. 
Let us prove that the latter is impossible.
Suppose that $\dim\, H^1(J(X),\Phi^2_3)=1$. Then, by (\ref{ex-3}),
\bea
&&
\frac{H^0(J(X),\Phi^2_3)}{H^0(J(X),d\gr_2\,\Omega^1)}\simeq H^0(J(X),R^2_3).
\label{8-3}
\ena
Since
\bea
&&
H^0(J(X),\Phi^2_3)=\ker\left(d:H^0(J(X),\gr_3\Omega^2)\lar H^0(J(X),\gr_4\,\Omega^3)\right),
\non
\ena
there exists an element $w$ of $H^0(J(X),\Omega^2(3))$ such that 
its image in $H^0(J(X),R^2_3)$ becomes a basis of $H^0(J(X),R^2_3)$. 
We assume $w=1/3\varphi^2_3$ in $(R^2_3)_0$. 
By Lemma \ref{lem-r2-r3}
$dw$ is a basis of $H^0(J(X),R^3_3)$ satisfying $dw=\varphi^3_3$ in $(R^3_3)_0$.
Let us write $dw=F_3\dum$ with $F_3\in A(3)$. 
Similarly to the proof of Lemma \ref{lem-41}, one can prove that
$F_3$ is contained in the space
\bea
&&
\sum\partial_i\partial_j A(2)+\sum_{i=1}^3\sum_{k=1}^5{\mathbb C}
v^k_i+M
\non
\\
&&
M:=\sum_{k=1}^5{\mathbb C}v^k+\sum_{i=1}^3 \partial_iA(2)+A(2).
\ena
Explicitly $F_3$ can be expressed as a linear combination of
\bea
&&
\zeta_{ijkl}, \zeta_{ijk}, \zeta_{ij},
 v^0_{ij}, v^0_i, v^0, v^m_i, v^m, 1,
\qquad
i,j,k,l\in\{1,2,3\},
\quad
1\leq m\leq 5.
\label{8-3}
\ena

Define the weight of $u_i$ to be $-(2i-1)$;
\bea
&&
\wt(u_i)=-(2i-1).
\non
\ena
In general we say that an element $a$ of $A$ has weight $d$ if
\bea
&&
a=\frac{a_{-6n+d}+(\hbox{lower weight terms})}{\sigma^n},
\quad
a_{-6n+d}\neq 0,
\quad
\wt\,a_i=i.
\non
\ena
By a calculation using (\ref{g3-sigma}) we have
\bea
&&
\wt\, \zeta_{i_1...i_k}=\sum_{j=1}^k(2i_j-1),
\non
\\
&&
\wt\,v^0=12,
\quad
\wt\,v^1=8,
\quad
\wt\,v^2=10,
\quad
\wt\,v^3=12,
\quad
\wt\,v^4=14,
\quad
\wt\,v^5=16.
\non
\ena
Since
\bea
F_3=\frac{1+(\hbox{lower weight terms})}{\sigma^3},
\non
\ena
we have $\wt\,F_3=18$.
Elements with weights no less than $18$ among (\ref{8-3}) are
\bea
&&
v^0_{33}(22),
\quad
v^5_3(21),
\quad
v^0_{23}(20),
\quad
\zeta_{3333}(20),
\quad
v^5_2(19),
\quad
v^4_3(19),
\quad
\zeta_{2333}(18),
\non
\\
&&
v^0_{13}(18),
\quad
v^0_{22}(18),
\label{83+}
\ena
where the number inside the bracket signifies the weight 
of the element.
A direct calculation shows
\bea
\sigma^4v^0_{33}&=&6u_1^2+\cdots,
\label{8-4------}
\\
\sigma^4v^5_{3}&=&-6u_1^3+\cdots,
\label{8-4-----}
\\
\sigma^4v^0_{23}&=&-2u_1^4-12u_1u_2+\cdots,
\label{8-4----}
\\
\sigma^4\zeta_{3333}&=&-6u_1^4+\cdots,
\label{8-4---}
\\
\sigma^4v^5_{2}&=&2u_1^5+12u_1^2u_2+\cdots,
\label{8-4--}
\\
\sigma^4v^4_{3}&=&-2u_1^5+6u_1^2u_2+\cdots,
\label{8-4-}
\\
\sigma^4\zeta_{2333}&=&2u_1^6+12u_1^3u_2+\cdots,
\label{8-4}
\\
\sigma^4v^0_{13}&=&\frac{4}{5}u_1^6-6u_1^3u_3+6u_1u_3-2S(u)+\cdots,
\label{8-5}
\\
\sigma^4v^0_{22}&=&\frac{2}{3}u_1^6+8u_1^3u_2+24u_2^2+4S(u)+\cdots,
\label{8-6},
\ena
where $\cdots$ signifies the lower weight terms.
Suppose that
\bea
&&
F_3=\sum_{i\leq j\leq k\leq l}c^1_{ijkl}\zeta_{ijkl}
+\sum_{i\leq j}c^2_{ij}v^0_{ij}+\sum_{i=1}^3\sum_{k=1}^5c^3_{k,i}v^k_i
\quad
\hbox{mod.} M.
\label{8-7}
\ena
Notice that 
\bea
&&
\sigma^4F_3=S(u)+\cdots,
\label{8-8}
\ena
and all elements in $M$ have weights less than $18$.
Multiply $\sigma^4$ to (\ref{8-7}) and compare homogeneous components from
the highest weight. 
Then the terms with weights greater than $-6=18-24$ should vanish
in the right hand side of (\ref{8-7}). It means that the following
coefficients are zero; 
\bea
&&
c^2_{33}, c^3_{5,3}, c^2_{23}, c^1_{3333}, c^3_{5,2}, c^3_{4,3}.
\non
\ena
For the weight $-6$ terms the following equation must hold:
\bea
S=c^1_{2333}P_1+c^2_{13}(P_2-2S)+c^2_{22}(P_3+4S),
\label{8-9}
\ena
where $P_1$, $P_2-2S$ and $P_3+4S$ are 
the top terms of the right hand side of 
(\ref{8-4}), (\ref{8-5}) and (\ref{8-6}) respectively.
However one can easily verify the linear independence of $S$, $P_1$,...,$P_3$.
Thus (\ref{8-9}) can not hold and consequently (\ref{8-2}) is proved. \qed
\vskip5mm

Recall that $(123;123)\in A(4)$ by Lemma \ref{ord-zeta} (iii) and
\bea
&&
d\zeta_1 \wedge d\zeta_2 \wedge d\zeta_3=(123;123) \dum.
\non
\ena

\begin{cor}\label{cor-81}
There exists a $2$-form $\xi\in H^0(J(X),\Omega^2(3))$ such that
$d\zeta_1 \wedge d\zeta_2 \wedge d\zeta_3-d\xi$ is contained
in $H^0(J(X),\Omega^3(3))$ and gives a basis of $H^0(J(X),R^3_3)$.
\end{cor}
\vskip2mm
\noindent
{\it Proof.}
Since we have already proved the equation (\ref{8-1}), for the first statement
 it is 
sufficient to prove 
\bea
&&
d\zeta_1 \wedge d\zeta_2 \wedge d\zeta_3\in H^0(J(X),d\gr_3\,\Omega^2).
\label{8-10}
\ena
Notice that 
\bea
d\zeta_1 \wedge d\zeta_2 \wedge d\zeta_3&=&
\frac{1}{3}d\tilde{\xi},
\non
\\
\tilde{\xi}&=&
\zeta_1d\zeta_2 \wedge d\zeta_3+\zeta_2d\zeta_3 \wedge d\zeta_1+
\zeta_3d\zeta_1 \wedge d\zeta_2.
\non
\ena
By a direct calculation we easily see that the order of poles of $\tilde{\xi}$
on $\Theta$ is at most $3$. Thus (\ref{8-10}) is proved.

Next let us prove that the coefficient of $\dum$ of 
$d\zeta_1 \wedge d\zeta_2 \wedge d\zeta_3-d\xi$
has a non-zero component of ${\mathbb C}F_3$ in the decomposition of 
Corollary \ref{cor-72}. Assume that this is not the case.
Then $(123;123)-\sum_{i=1}^3 \partial_i\xi_i$ can be written as a linear
combination of $\wp_{ij}$, $(ij;kl)$'s modulo 
$\sum_{i=1}^3\partial_i A$, where
\bea
&&
\xi=\xi_1 du_2 \wedge du_3+\xi_2 du_3 \wedge du_1+\xi_3 du_1 \wedge du_2.
\label{2form-xi}
\ena
It contradicts the linear independence of the basis of 
$A/\sum_{i=1}^3\partial_i A$ given by Proposition \ref{cohomology-basis}.
\qed

\begin{cor}\label{cor-82}
The ${\cal D}$-module is generated by $A(3)$.
\end{cor}
\vskip2mm
\noindent
{\it Proof.}
Let $F_4\dum$, $F_4\in A(4)$ be a basis of $H^0(J(X),R^3_4)$.
Then, using Proposition \ref{prop-81}, in a similar way to the proof of
Proposition \ref{gen-4} one can prove that $F_4$ is contained in 
${\cal D}A(3)$. Thus $A(4)$ is contained in ${\cal D}A(3)$. By Proposition
\ref{gen-4} we have $A={\cal D}A(3)$.
\qed

\begin{cor}\label{cor-83}
The ${\cal D}$-module $A$ is generated by 
representatives (\ref{c-basis}) of 
$A/\sum_{i=1}^3\partial_i A$. Explicitly
\bea
&&
A={\cal D}1+\sum_{i\leq j}{\cal D}\zeta_{ij}+
\sum{\cal D}(ij;kl)
+{\cal D}(123;123),
\label{8-11}
\ena
where the sum for $(ij;kl)$ is over elements appeared in 
(\ref{c-basis}).
\end{cor}
\vskip2mm
\noindent
{\it Proof.}
By Corollary \ref{basis-a2}, \ref{cor-72} and \ref{cor-81}, it is sufficient
to prove that
\bea
&&
(123;123)-\sum_{i=1}^3\partial_i\xi_i\in 
\hbox{the right hand side of (\ref{8-11}),}
\non
\ena
for a two form $\xi$ as in Corollary \ref{cor-81}.
To this end it is sufficient to prove
\bea
&&
\xi_i\in {\mathbb C}+\sum_{i\leq j}{\mathbb C}\zeta_{ij}+
\sum_{i=1}^3{\mathbb C}v^0_i+\sum_{k=1}^5{\mathbb C}v^k,
\label{8-11+}
\ena
for $1\leq i\leq 3$.
Let $F_3\dum$, $F_3\in A(3)$ be a basis of $H^0(J(X),R^3_3)$
such that
\bea
&&
F_3=\frac{1+(\hbox{lower degree terms})}{\sigma^3}.
\non
\ena

Let us write $\xi_i$ as a linear combination of a basis of $A(3)$:
\bea
&&
\xi_i=c^1_iF_3+\sum_{k=1}^5 c^2_k v^k+\sum_{k=1}^3 c^3_k v^0_k+
\sum_{j\leq k\leq l}c^4_{jkl}\zeta_{jkl}+
\sum_{j\leq k}c^5_{jk}\zeta_{jk}
+c^6v^0+c_7.
\non
\ena
Then
\bea
&&
\partial_i\xi_i=c^1_i\partial_iF_3+\sum_{k=1}^5 c^2_k v^k_i
+\sum_{k=1}^3 c^3_k v^0_{ik}+
\sum_{j\leq k\leq l}c^4_{jkl}\zeta_{ijkl}+
\sum_{j\leq k}c^5_{jk}\zeta_{ijk}
+c^6v^0_i.
\non
\ena
Among elements appeared in the right hand side of this equation,
those with the weights no less than $18$ are (\ref{83+}) and
\bea
&&
\partial_1 F_3 (19),
\quad
\partial_2 F_3 (21),
\quad
\partial_3 F_3  (23).
\non
\ena
By a calculation we have
\bea
\sigma^4\partial_3 F_3&=&-3u_1+\cdots,
\label{8-12}
\\
\sigma^4\partial_2 F_3&=&u_1^3+6u_2+\cdots,
\label{8-13}
\\
\sigma^4\partial_1 F_3&=&-\frac{2}{5}u_1^5+3u_1^2u_2-3u_3+\cdots.
\label{8-14}
\ena
Since $(123;123)-\sum_{i=1}^3 \partial_i\xi_i\in A(3)$, the weights 
of elements of $A(3)$ are at most $18$ and $\wt\,(123;123)=18$,
we have 
\bea
&&
\left(\sigma^4\sum_{i=1}^3\partial_i\xi_i\right)_{\geq -5}=0,
\non
\ena
where $(\quad)_{\geq -5}$ denotes the terms with weights 
no less than than $-5$.
Using the expansion (\ref{8-12})-(\ref{8-14}) and (\ref{8-4------})-(\ref{8-6}),we easily find $c^1_i=0$. Thus (\ref{8-11+}) is proved.
\qed

\section{Linear Basis of Abelian Functions}
In this section we determine a basis of $A$ as a vector space. 
It is constructed as a subset of the set of derivatives of the basis
of $A/\sum_{i=1}^3 \partial_i A$ given in Proposition \ref{cohomology-basis}.

The obvious relation $d(d\zeta_i\wedge d\zeta_j)=0$ gives
\bea
&&
\partial_3(ij;12)-\partial_2(ij;13)+\partial_1(ij;23)=0.
\non
\ena
Thus we have
\bea
&&
\partial_3(12;ij)=\partial_2(13;ij)-\partial_1(23;ij),
\label{3-der}
\ena
for any $i<j$.
Using these relations it is possible to erase $u_3$-derivatives from the
derivatives of $(12;ij)$.

We use the notation $(1^{a_1}2^{a_2}3^{a_3})=(1,...,1,2,...,2,3,...,3)$ where
$i$ appears $a_i$ times, and
\bea
&&
\zeta_{1^{a_1}2^{a_2}3^{a_3}}=
\partial_1^{a_1}\partial_2^{a_2}\partial_3^{a_3}\log\,\sigma,
\quad
w_{1^{a_1}2^{a_2}3^{a_3}}=\partial_1^{a_1}\partial_2^{a_2}\partial_3^{a_3}w
\non
\ena
for $w\in A$.

\begin{theorem}\label{lin-basis}
The following elements give a basis of $A$ as a vector space:
\bea
&&
1,\quad
\zeta_{1^{a_1}2^{a_2}3^{a_3}}
\quad
(a_1+a_2+a_3\geq 2),
\non
\\
&&
(12;12)_{1^{a_1}2^{a_2}},
\quad
(12;13)_{1^{a_1}2^{a_2}},
\quad
(12;23)_{1^{a_1}2^{a_2}},
\quad
(a_1,a_2\geq 0),
\non
\\
&&
(13;13)_{1^{a_1}2^{a_2}3^{a_3}},
\quad
(13;23)_{1^{a_1}2^{a_2}3^{a_3}},
\quad
(23;23)_{1^{a_1}2^{a_2}3^{a_3}}
\quad
(a_1, a_2, a_3 \geq 0),
\non
\\
&&
(123;123)_{1^{a_1}2^{a_2}3^{a_3}}
\quad
(a_1, a_2, a_3 \geq 0).
\label{c-basis1}
\ena
\end{theorem}
\vskip2mm
\noindent
The elements (\ref{c-basis1}) generate 
$A$ as a vector space by Corollary \ref{cor-83} and relations (\ref{3-der}). 
Therefore we have to prove the linear independence 
of the elements (\ref{c-basis1}) in order to prove the theorem.

Let $\xi$ be an element of $H^0(J(X),\Omega^2(3))$ as in Corollary \ref{cor-81}
and $\xi_i$ its components defined by (\ref{2form-xi}). We set
\bea
&&
u^3=(123;123)-\sum_{i=1}^3\partial_i\xi^i.
\non
\ena

Then $u^3\in A(3)$, $u^3\dum$ is a basis of $H^0(J(X),R^3_3)$ and $\partial_i\xi^i$ 
is a linear combination of elements in (\ref{c-basis1}) other than 
$(123;123)_{1^{a_1}2^{a_2}3^{a_3}}$'s by (\ref{8-11+}).

Thus the theorem is equivalent to saying that the following
elements are a ${\mathbb C}$-basis of $A$:
\bea
&&
1,\quad
\zeta_{1^{a_1}2^{a_2}3^{a_3}}
\quad
(a_1+a_2+a_3\geq 2),
\non
\\
&&
(12;12)_{1^{a_1}2^{a_2}},
\quad
(12;13)_{1^{a_1}2^{a_2}},
\quad
(12;23)_{1^{a_1}2^{a_2}},
\quad
(a_1,a_2\geq 0),
\non
\\
&&
(13;13)_{1^{a_1}2^{a_2}3^{a_3}},
\quad
(13;23)_{1^{a_1}2^{a_2}3^{a_3}},
\quad
(23;23)_{1^{a_1}2^{a_2}3^{a_3}}
\quad
(a_1, a_2, a_3 \geq 0),
\non
\\
&&
u^3_{1^{a_1}2^{a_2}3^{a_3}}\, (a_1, a_2, a_3 \geq 0).
\label{c-basis2}
\ena
We shall prove the linear independence of them. To this end we prepare
some lemmas.

Let us define $u^n\in A(n)$, $n\geq 4$ satisfying 
${\mathbb C}u^n\dum=H^0(J(X),R^3_n)$ inductively as follows.

For each $n\geq 2$ let $B_n$ be the set of elements
\bea
&&
\zeta_{1^{a_1}2^{a_2}3^{a_3}}
\quad
(a_1+a_2+a_3=n),
\non
\\
&&
\left((13;13)-(12;23)\right)_{1^{a_1}2^{a_2}3^{a_3}},
\quad
(a_1+a_2+a_3=n-2),
\non
\\
&&
(12;12)_{1^{a_1}2^{a_2}},
\quad
(12;13)_{1^{a_1}2^{a_2}},
\quad
(12;23)_{1^{a_1}2^{a_2}},
\quad
(a_1+a_2=n-3),
\non
\\
&&
(13;23)_{1^{a_1}2^{a_2}3^{a_3}},
\quad
(23;23)_{1^{a_1}2^{a_2}3^{a_3}}
\quad
(a_1+a_2+a_3=n-3),
\non
\\
&&
u^3_{1^{a_1}2^{a_2}3^{a_3}}\, (a_1+a_2+a_3=n-3).
\non
\ena
For $n=2$ we understand that there are no elements specified by the condition
$a_1+a_2+a_3=n-3$ or $a_1+a_2=n-3$. We set $B_0=\{1\}$, 
$B_1=\emptyset$.
Notice that $B_n$ is a basis of $\gr_n\,A$ for $n\leq 3$ 
by Corollary \ref{basis-a2} and \ref{cor-72}.

By Lemma \ref{lem-41}, Corollary \ref{cor-82} there is a linear combination
$P^4$ of elements in $B_5$ such that $u^4:=P^4$ is an element of $A(4)$ and 
$u^4\dum$ is a basis of $H^0(J(X),R^3_4)$.
Suppose that $u^i$, $i\leq n$ are defined. 
By Lemma \ref{lem-41} there is a linear combination $P^{n+1}$ of elements
in $B_{n+2}\cup\{u^j_{1^{a_1}2^{a_2}3^{a_3}}| 4\leq j\leq n, 
a_1+a_2+a_3=n+2-j\}$ such that $u^{n+1}:=P^{n+1}$
is an element of $A(n+1)$ and $u^{n+1}\dum$ gives a basis of $H^0(J(X),R^3_{n+1})$.

\begin{lemma}\label{lem-91} For $n\geq 4$,  
$\{u^j_{1^{a_1}2^{a_2}3^{a_3}}|\,4\leq j\leq n,a_1+a_2+a_3=n-j\}$
is linearly independent as elements of $A$.
\end{lemma}

The proof of the lemma is given later.

Let $P^j_{1^{a_1}2^{a_2}3^{a_3}}$ denote 
$\partial_{1^{a_1}}\partial_{2^{a_2}}\partial_{3^{a_3}}P^j$, 
where $u_3$-derivatives of $(12;ij)$'s are erased by the relation
(\ref{3-der}). Therefore $P^j_{1^{a_1}2^{a_2}3^{a_3}}$, $a_1+a_2+a_3=n-j$
 is a linear combination
of elements in $B_{n+1}\cup\{u^j_{1^{a_1}2^{a_2}3^{a_3}}| 4\leq j\leq n+1, 
a_1+a_2+a_3=n+1-j\}$. Let
\bea
&&
C_n=B_n\sqcup \{u^j_{1^{a_1}2^{a_2}3^{a_3}}|\, 4\leq j\leq n, 
a_1+a_2+a_3=n-j\}.
\non
\ena
Let us consider the set of symbols 
\bea
&&
\bar{C}_n=\bar{B}_n\sqcup \{\bar{u}^j_{1^{a_1}2^{a_2}3^{a_3}}|\, 4\leq j\leq n, a_1+a_2+a_3=n-j\},
\non
\ena
where $\bar{B}_n$ is the set of symbols obtained by making a bar to each element of $B_n$, $\bar{B}_n=\{\bar{\zeta}_{1^{a_1}2^{a_2}3^{a_3}},...\}$. 
The elements in 
$\bar{C}_n$ are considered to be linearly independent. For a linear combination
$P$ of elements in $C_n$, let $\bar{P}$ be the linear combination of elements
in $\bar{C}_n$ obtained by making a bar to each element of $C_n$ 
appeared in $P$. A priori $\bar{P}$ may not be uniquely defined since the
expression $P$ may not be unique. Take any one of the expression and
make $\bar{P}$. We denote the vector space with the elements of $\bar{C}_n$ as a basis by $\hbox{Span}_{\mathbb{C}}\bar{C}_n$.

\begin{lemma}\label{lem-92}
The set $\{\bar{P}^j_{1^{a_1}2^{a_2}3^{a_3}}|\,4\leq j\leq n,\,
a_1+a_2+a_3=n-j\}$ is linearly independent in 
$\hbox{Span}_{\mathbb{C}}\bar{C}_{n+1}$ and
\bea
&&
\gr^P_{n+1}\,A\simeq \frac{\hbox{Span}_{\mathbb{C}}\bar{C}_{n+1}}
{\oplus_{4\leq j\leq n, a_1+a_2+a_3=n-j} 
{\mathbb C}\bar{P}^j_{1^{a_1}2^{a_2}3^{a_3}}},
\label{gr-(n+1)}
\ena
where the map from the RHS to the LHS is defined simply by erasing 
bars of symbols.
\end{lemma}

Lemma \ref{lem-92} follows from Lemma \ref{lem-91}.
In fact assume Lemma \ref{lem-91}. Then
\bea
\sum_{4\leq j\leq n,a_1+a_2+a_3=n-j}
 c^{a_1a_2a_3}_j\bar{P}^j_{1^{a_1}2^{a_2}3^{a_3}}=0,
&\Rightarrow&
\sum_{}
 c^{a_1a_2a_3}_jP^j_{1^{a_1}2^{a_2}3^{a_3}}=0,
\non
\\
&\Rightarrow&
\sum c^{a_1a_2a_3}_ju^j_{1^{a_1}2^{a_2}3^{a_3}}=0,
\non
\\
&\Rightarrow&
c^{a_1a_2a_3}_j=0,
\non
\ena
since $\{u^j_{1^{a_1}2^{a_2}3^{a_3}}|4\leq j\leq n, a_1+a_2+a_3=n-j\}$ 
is linearly independent by Lemma \ref{lem-91}.
Thus $\{\bar{P}^j_{1^{a_1}2^{a_2}3^{a_3}}|4\leq j\leq n, a_1+a_2+a_3=n-j\}$
 is linearly independent.

Next let us prove (\ref{gr-(n+1)}).
We already know that the map given there is well defined, surjective and
$\dim\,\gr_{n+1}\,A=(n+1)^3-n^3$.
Let us compute the dimension of the RHS. We denote by ${}_nH_r$ the number of
combinations taking $r$ elements from $n$ elements admitting repetition.
Then
\bea
&&
\sharp \bar{C}_n={}_3H_{n}+{}_3H_{n-2}+3{}_2H_{n-3}+3{}_3H_{n-3}+
\sum_{j=4}^n{}_3H_{n-j},
\non
\\
&&
\sharp\{\bar{P}^j_{1^{a_1}2^{a_2}3^{a_3}}|\,
4 \leq j\leq n-1, a_1+a_2+a_3=n-1-j\}=\sum_{j=4}^{n-1}{}_3H_{n-1-j},
\non
\ena
and
\bea
\dim (\hbox{RHS of (\ref{gr-(n+1)})})&=&
\sharp \bar{C}_{n+1}
-
\sharp\{\bar{P}^j_{1^{a_1}2^{a_2}3^{a_3}}|\,
4 \leq j\leq n, a_1+a_2+a_3=n-j\}
\non
\\
&=&
{}_3H_{n+1}+{}_3H_{n-1}+3{}_2H_{n-2}+3{}_3H_{n-2}+{}_3H_{n-3}
\non
\\
&=&
(n+1)^3-n^3.
\non
\ena
Thus the map is an isomorphism.
\qed

\vskip5mm
\noindent
Proof of Lemma \ref{lem-91}
\vskip2mm
\noindent
Let us prove the lemma by the induction on $n$.
Consider the case of $n=4$. 
In this case 
$\{u^j_{1^{a_1}2^{a_2}3^{a_3}}|\,4\leq j\leq n, a_1+a_2+a_3=n-j\,\}=\{u^4\}$. By the definition 
$u^4\neq 0$ and the lemma is obvious.

Assume the lemma until $n$. By the induction hypothesis for $n$ we have the
isomorphism (\ref{gr-(n+1)}). Then
\bea
&&\sum_{4\leq j\leq n+1,a_1+a_2+a_3=n+1-j}
 c^{a_1a_2a_3}_ju^j_{1^{a_1}2^{a_2}3^{a_3}}=0\quad\hbox{in $A$},
\non
\\
&&
\qquad\qquad
\Rightarrow
\sum c^{a_1a_2a_3}_ju^j_{1^{a_1}2^{a_2}3^{a_3}}=0\quad\hbox{in $\gr_{n+1}\,A$},
\non
\\
&&
\qquad\qquad
\Rightarrow
\sum c^{a_1a_2a_3}_j\bar{u}^j_{1^{a_1}2^{a_2}3^{a_3}}=0\quad
\hbox{in RHS of (\ref{gr-(n+1)})}.
\non
\ena
The last equation implies
\bea
&&
\sum c^{a_1a_2a_3}_j\bar{u}^j_{1^{a_1}2^{a_2}3^{a_3}}=
\sum_{4\leq j\leq n,a_1+a_2+a_3=n-j} \tilde{c}^{a_1a_2a_3}_j
\bar{P}^j_{1^{a_1}2^{a_2}3^{a_3}}
\quad\hbox{in $\bar{C}_{n+1}$},
\label{rel-cbar}
\ena
for some constants $\tilde{c}^{a_1a_2a_3}_j$. Then we have
\bea
&&
0=\sum c^{a_1a_2a_3}_j u^j_{1^{a_1}2^{a_2}3^{a_3}}=
\sum_{4\leq j\leq n,a_1+a_2+a_3=n-j}
 \tilde{c}^{a_1a_2a_3}_j u^j_{1^{a_1}2^{a_2}3^{a_3}}
\quad\hbox{in $A$},
\non
\ena
which implies
\bea
&&
\tilde{c}^{a_1a_2a_3}_j=0,
\non
\ena
since $\{u^j_{1^{a_1}2^{a_2}3^{a_3}}\,|\,4\leq j\leq n,a_1+a_2+a_3=n-j\}$ 
is linearly independent by the hypothesis
of induction. By (\ref{rel-cbar}) we have
\bea
&&
\sum_{4\leq j\leq n+1,a_1+a_2+a_3=n+1-j}
 c^{a_1a_2a_3}_j \bar{u}^j_{1^{a_1}2^{a_2}3^{a_3}}=0.
\non
\ena
Thus all $c^{a_1a_2a_3}_j=0$ and 
$\{u^j_{1^{a_1}2^{a_2}3^{a_3}}|\,4\leq j\leq n+1,a_1+a_2+a_3=n+1-j\}$ 
is linearly independent. \qed

\vskip5mm
\noindent
Proof of Theorem \ref{lin-basis}
\vskip2mm
\noindent
The linear independence of elements (\ref{c-basis2}) is equivalent to
that of elements of $B:=\sqcup_{n=0}^\infty B_n$. We prove that
elements of $B$ are linearly independent.

Consider the linear relation among elements of $B$ and write it as
\bea
&&
Q_n+Q_{n-1}+\cdots+Q_0=0,
\label{eq-90}
\ena
where $Q_i$ is a linear combination of elements in $B_i$.
Then
\bea
&&
Q_n=0\quad\hbox{in $\gr_n\,A$},
\non
\ena
and 
\bea
&&
\bar{Q}_n=\sum_{4\leq j\leq n-1,a_1+a_2+a_3=n-1-j}
 c^{n;a_1a_2a_3}_j \bar{P}^j_{1^{a_1}2^{a_2}3^{a_3}}
\quad
\hbox{in ${\bar C}_n$},
\non
\ena
for some constants $c^{n;a_1a_2a_3}_j$. It implies
\bea
&&
Q_n=
\sum_{4\leq j\leq n-1,a_1+a_2+a_3=n-1-j}
 c^{n;a_1a_2a_3}_j u^j_{1^{a_1}2^{a_2}3^{a_3}}
\quad
\hbox{in $A$}.
\label{n10-1}
\ena
Notice that the right hand side of (\ref{n10-1})
 is a linear combination of elements in $C_{n-1}$.
We have
\bea
&&
\left(\sum_{4\leq j\leq n-1,a_1+a_2+a_3=n-1-j} 
c^{n;a_1a_2a_3}_j u^j_{1^{a_1}2^{a_2}3^{a_3}}+Q_{n-1}\right)
+Q_{n-2}+\cdots+Q_0=0.
\non
\ena
Similarly, for $k\leq n$, there are constants 
$c^{k;a_1a_2a_3}_j$, $4\leq j\leq k-1$,
 $a_1+a_2+a_3=k-1-j$, such that
\bea
&&
\sum  c^{k;a_1a_2a_3}_j\bar{u}^j_{1^{a_1}2^{a_2}3^{a_3}}+\bar{Q}_{k-1}
=\sum c^{k-1;a_1a_2a_3}_j \bar{P}^j_{1^{a_1}2^{a_2}3^{a_3}}
\quad\hbox{in $\bar{C}_{k-1}$},
\label{eq-91}
\\
&&
\sum c^{k-1;a_1a_2a_3}_j u^j_{1^{a_1}2^{a_2}3^{a_3}}+Q_{k-2}+\cdots+Q_0=0
\quad
\hbox{in $A$.}
\label{eq-92}
\ena
Taking $k=6$ in (\ref{eq-92}) we have
\bea
&&
c^{5;000}_{4}u^4+Q_4+\cdots+Q_0=0.
\non
\ena
It implies
\bea
&&
c^{5;000}_{4}=0,
\quad
\bar{Q}_i=0, 
\quad 
0\leq i\leq 4,
\non
\ena
since $\sqcup_{i=0}^4 B_i\sqcup \{u^4\}$ is linearly independent.
Then, by (\ref{eq-91}),
\bea
&&
\sum c^{6;a_1a_2a_3}_j \bar{u}^j_{1^{a_1}2^{a_2}3^{a_3}}+\bar{Q}_{5}
=0,
\non
\ena
which implies 
\bea
&&
c^{6;a_1a_2a_3}_j=0,
\quad
\bar{Q}_5=0,
\non
\ena
since 
$\{\bar{u}^j_{1^{a_1}2^{a_2}3^{a_3}}\,|4\leq j\leq 5, a_1+a_2+a_3=5-j\,\}\cap \bar{B}_5=\emptyset$ 
in $\bar{C}_5$.
Repeating similar arguments we have $\bar{Q}_i=0$ for any $i$. It means that 
the linear relation (\ref{eq-90}) is trivial. Thus $B$ is linearly
independent. \qed

\section{Proof of Theorem \ref{main-th}}
\vskip2mm
\noindent
By Lemma \ref{lem-n22} and (\ref{eq-n25}) we have
\bea
&&
(k_1,...,k_m;l_1,...,l_m)_{1^{a_1}2^{a_2}3^{a_3}}\in A_d,
\quad
d=2\sum_{i=1}^m(k_i+l_i-1)+a_1+3a_2+5a_3.
\non
\ena
Let $B^n$, $n\geq 1$, be the subset of (\ref{c-basis1}) consisting of 
elements of the form
\bea
&&
(k_1,...,k_m;l_1,...,l_m)_{1^{a_1}2^{a_2}3^{a_3}},
\quad
2\sum_{i=1}^m(k_i+l_i-1)+a_1+3a_2+5a_3=n.
\non
\ena
We set $B^0=\{1\}$. For example
\bea
&&
B^1=\emptyset,
\quad
B^2=\{\zeta_{11}\},
\quad
B^3=\{\zeta_{1^3}\},
\quad
B^4=\{\zeta_{1^4},\zeta_{12}\},
\non
\\
&&
B^5=\{\zeta_{1^5},\zeta_{112}\}
\quad
B^6=\{\zeta_{1^6},\zeta_{1^32},\zeta_{13},\zeta_{22}\},
\quad
B^7=\{\zeta_{1^7},\zeta_{1^42},\zeta_{113},\zeta_{122}\},
\non
\\
&&
B^8=\{\zeta_{1^8},\zeta_{1^52},\zeta_{1^33},\zeta_{1^22^2},\zeta_{23},(12;12)\}.\non
\ena

\begin{lemma}\label{lem-11-1} We have
\bea
&&
\sum_{n=0}^\infty q^n\dim\,B^n=\ch(\gr^{KP}\,A).
\non
\ena
\end{lemma}
\vskip2mm
\noindent
{\it Proof.}
We have
\bea
\sum_{n=0}^\infty q^n\dim\,B^n
&=&
1+\sum_{a_1+a_2+a_3\geq 2}q^{a_1+3a_2+5a_3}
+(q^8+q^{10}+q^{12})\sum_{a_1,a_2=0}^\infty q^{a_1+3a_2}
\non
\\
&\quad&
+(q^{12}+q^{14}+q^{16}+q^{18})\sum_{a_1,a_2,a_3=0}^\infty q^{a_1+3a_2+5a_3}
\non
\\
&=&
-q-q^3-q^5+
\frac{1+q^{12}+q^{14}+q^{16}+q^{18}}{(1-q)(1-q^3)(1-q^5)}
+\frac{q^8+q^{10}+q^{12}}{(1-q)(1-q^3)}.
\label{11-2}
\ena
On the other hand we have 
\bea
&&
\ch(\gr^{KP}\,A)=
\frac{[\,\frac{1}{2}\,]_{q^2}[7]_{q^2}!}
{[3]_{q^2}![4]_{q^2}![3+\frac{1}{2}]_{q^2}!}.
\label{11-3}
\ena
by Theorem \ref{th-n24}. By a direct calculation one can show
that (\ref{11-2}) and (\ref{11-3}) are equal.
\qed

\begin{lemma}\label{lem-11-2}
The set $B^n$ is a basis of $\gr^{KP}_n\,A$.
\end{lemma}
\vskip2mm
\noindent
{\it Proof.}
The lemma can be easily proved by induction on $n$ 
using the linear independence of (\ref{c-basis1}) and
Lemma \ref{lem-11-1}.
\qed
\vskip5mm

It follows from this lemma that $\gr^{KP}\,A$ is generated
by $1$, $(ij;kl)$, $(123;123)$ over ${\cal D}$. Thus
Theorem \ref{main-th} is proved. \qed

\begin{cor}
The following set of elements is a basis of $\gr_n^{KP}\,A$ for $n\geq 2$;
\bea
&&
\zeta_{1^{a_1}2^{a_2}3^{a_3}}
\quad
(a_1+3a_2+5a_3=n),
\quad
(12;12)_{1^{a_1}2^{a_2}}
\quad
(8+a_1+3a_2=n),
\non
\\
&&
(12;13)_{1^{a_1}2^{a_2}}
\quad
(10+a_1+3a_2=n),
\non
\\
&&
(12;23)_{1^{a_1}2^{a_2}}
\quad
(12+a_1+3a_2=n),
\non
\\
&&
(13;13)_{1^{a_1}2^{a_2}3^{a_3}}
\quad
(12+a_1+3a_2+5a_3=n),
\non
\\
&&
(13;23)_{1^{a_1}2^{a_2}3^{a_3}}
\quad
(14+a_1+3a_2+5a_3=n),
\non
\\
&&
(23;23)_{1^{a_1}2^{a_2}3^{a_3}}
\quad
(16+a_1+3a_2+5a_3=n),
\non
\\
&&
(123;123)_{1^{a_1}2^{a_2}3^{a_3}}
\quad
(18+a_1+3a_2+5a_3=n).
\non
\ena
\end{cor}

\section{Concluding Remarks}
In this paper we have determined the ${\cal D}$-module structure of the
affine ring $A$ of the affine Jacobian of a hyperelliptic curve of genus
$3$. The ${\cal D}$-free resolution conjectured in the paper \cite{NS1} 
is proved to be true. In particular generators and relations among them 
over ${\cal D}$ are determined. 
A ${\mathbb C}$-linear basis of $A$ is also given in terms
of Klein's hyperelliptic $\wp$-functions.

Two filtrations, pole and KP filtrations, are introduced for $A$.
It is proved that the graded ring $\gr^P\,A$ associated with 
the pole filtration is not finitely generated. 
The reason is the existence of the singularity of the 
theta divisor. We study the effect of the singularity in detail and
reveal the structure on how $A$ becomes finitely generated
although $\gr^P\,A$ is not. We think it a typical structure valid for
hyperelliptic Jacobians of genus $g\geq 3$ or more generally 
principally polarized abelian varieties with the singular theta divisors.
Unfortunately we could not find explicit formulae
for a ${\mathbb C}$-linear basis of $\gr^P_n\,A$ for each $n$.
It is an attractive problem to find them. The result will have
 an application to addition formulae of Frobenius-Stickelberger type.

The KP-filtration fits more naturally to the description of $A$. In fact
the graded ring $\gr^{KP}\,A$ associated with the KP-filtration is proved 
to be finitely generated and a ${\mathbb C}$-linear
basis of $\gr^{KP}\,A$ is explicitly constructed. In general KP-filtration
will be appropriate to describe $A$ of the affine Jacobian. However it
is effective to use both filtrations to prove something.

It is worth pointing that $\gr^{KP}\,A$ is isomorphic to
the ring $A_0$ corresponding to the most degenerate case,
that is, all coefficients $\lambda_i$
of the hyperelliptic curve are equal to zero \cite{NS1}. 
The latter ring is generated by logarithmic derivatives of 
a Schur function and thereby is a subring of the ring of 
rational functions. So Conjecture \ref{conj1} is reduced to
the problem on rational functions. We still do not know whether
it helps to prove the conjecture but expect it does.

Finally we remark that it is interesting to consider the deformation 
of the present case.
Namely consider the space of meromorphic sections of a non-trivial 
flat line bundle and determine the ${\cal D}$-module structure of it. 
The generic case had been studied in \cite{N0} in relation with
the problem of constructing commuting differential operators. 
It is curious to study whether the affine ring of an affine Jacobian can be
embedded in the ring of differential operators.

\vskip5mm
\noindent
{\large {\bf Acknowledgement}} 
\vskip3mm
\noindent
This research is supported by Grant-in-Aid for Scientific Research (B) 
17340048.
A part of the results of this paper was presented 
in "Integrable Systems, Geometry and Abelian Functions"
 held at Tokyo Metropolitan University, May 2005,
a series of lectures at University of Tokyo, November 2005
and "Methods of Integrable Systems in Geometry" held at 
University of Durham, August 2006. I would like to
thank Martin Guest, Yoshihiro \^Onishi, Shigeki Matsutani
 and Michio Jimbo for inviting me to give talks 
at Symposiums and lectures. 
I am also grateful to Koji Cho for fruitful discussions.


\begin{thebibliography}{9999}

\bibitem{A1}
M. Artin, On the solutions of analytic equations,
{\it Invent. Math.} {\bf 62} (1968), 277-291.

\bibitem{A2}
M. Artin, Algebraic approximation of structures over complete
local rings,
{\it Publ. IHES} {\bf 36} (1969), 23-58.

\bibitem{AH}
W. V. D. Hodge and M. F. Atiyah, Integrals of the second kind on an algebraic
variety,
{\it Annals of Math.} {\bf 62} (1955), 56-91.

\bibitem{B1}
H. F. Baker,
{\it Abelian functions},
1897, Cambridge University Press.

\bibitem{B2}
H. F. Baker,
On the hyperelliptic sigma functions,
{\it Amer. J. Math.} {\bf 20-4} (1898), 301-384.

\bibitem{B3}
H. F. Baker,
On a system of differential equations leading to periodic functions,
{\it Acta Math.} {\bf 27} (1903), 135-156.

\bibitem{BL}
V. M. Buchstaber and D. V. Leykin,
Addition laws on Jacobian varieties of plane algebraic curves,
{\it Proc. Steklov Inst. of Math. } {\bf 251} (2005), 1-72.

\bibitem{BEL1}
V. M. Buchstaber, V. Z. Enolskii and D. V. Leykin,
Kleinian functions, hyperelliptic Jacobians and
applications, in {\it Reviews in Math. and Math. Phys.}
{\bf Vol.10}, No.2, Gordon and Breach, London, 1997, 1-125.

\bibitem{CN}
K. Cho and A. Nakayashiki,
Differential structure of abelian functions,
{\it Int. J. Math. } {\bf 19} (2008), 145-171.


\bibitem{DJKM}
E. Date, M. Jimbo, M. Kashiwara and T. Miwa,
Transformation groups for soliton equations,
Nonlinear Integrable Systems-Classical Theory and Quantum Theory,
M. Jimbo and T. Miwa (eds.), World Scientific, Singapore, 1983,
39-119.

\bibitem{EG}
J. C. Eilbeck, V. Z. Enolskii and J. Gibbons,
Sigma, tau and abelian functions of algebraic curves,
J. Phys. A: Math. Theor. {\bf 43} (2010), 455216 (20pp)

\bibitem{EEMOP}
J. C. Eilbeck, V.Z. Enolskii, S. Matsutani, Y. \^Onishi and E. Previato,
Abelian functions for trigonal curves of genus three,
{\it Int. Math. Res. Not.} {\bf 2007} Article ID rnm140, 38pages (2007).

\bibitem{EEP}
J. C. Eilbeck, V. Z. Enolskii and E. Previato,
On a generalized Frobenius-Stickelberger addition formula,
{\it Lett. Math. Phys. } {\bf 65} (2003), 5-17.

\bibitem{K1}
F. Klein, Ueber hyperelliptische Sigmafunctionen,
{\it Math. Ann.} {\bf 27} (1886), 341-464.

\bibitem{K2}
F. Klein, Ueber hyperelliptische Sigmafunctionen 
(Zweiter Aufsatz),
{\it Math. Ann.} {\bf 32} (1888), 351-380.

\bibitem{M0}
D. Mumford,
{\it Abelian Varieties},
Oxford University Press, 1970.

\bibitem{M2}
D. Mumford,
{\it Tata lectures on theta} {\bf II},
Birkhauser, 1983.

\bibitem{N0}
A. Nakayashiki, Structure of Baker-Akhiezer modules of principally polarized
abelian varieties, commuting partial differential operators and associated 
integrable systems, {\it Duke Math. J.}{\bf 62-2} (1991) 315 --358.

\bibitem{N1}
A. Nakayashiki, On the cohomologies of theta divisors of 
hyperelliptic Jacobians, {\it Contemporary Math.} {\bf 309} (2002), 177-183.

\bibitem{N2}
A. Nakayashiki, Algebraic expressions of sigma functions for (n,s) curves, 
 Asian J. Math. {\bf 14-2} (2010), 175-212.

\bibitem{N3}
A. Nakayashiki, Sigma function as a tau function, 
Int. Math. Res. Not. {\bf 2010-3} (2010), 373-394.

\bibitem{NS1}
A. Nakayashiki and F. Smirnov, 
Cohomologies of affine hyperelliptic Jacobi varieties and integrable systems,
{\it Comm. Math. Phys.} {\bf 217} (2001), 623-652.

\bibitem{NS2}
A. Nakayashiki and F. Smirnov, 
Euler characteristics of theta divisors of Jacobians for spectral curves, 
{\it CRM Proc. and Lect. Notes} {\bf 32}, Vadim B. Kuznetsov, ed.
(2002), 239-246,


\bibitem{Sat}
M. Sato and Y. Sato, Soliton equations as dynamical systems on infinite
dimensional Grassmann manifold, Nolinear Partial Differential Equations
in Applied Sciences, P.D. Lax, H. Fujita and G. Strang (eds.), 
North-Holland, Amsterdam, and Kinokuniya, Tokyo, 1982, 259-271.

\bibitem{Z}
O. Zariski,
{\it Algebraic Surfaces}, Classic in Mathematics,
Springer-Verlag Berlin Heidelberg, 1995.


\end{thebibliography}
\end{document}